\documentclass[11pt]{amsart}
\usepackage{asymptote}
\usepackage{latexsym,amssymb,amsmath}
\usepackage{mathrsfs}
\usepackage{epsfig}
\usepackage{amsmath}
\usepackage{amsfonts}
\usepackage{amssymb}
\usepackage{amssymb,amscd}
\usepackage{graphicx}
\usepackage{subfigure}
\usepackage{psfrag}
\RequirePackage{graphics}

\begin{document}
\newtheorem{thm}{Theorem}
\newtheorem{lem}[thm]{Lemma}
\newtheorem{cor}[thm]{Corollary}
\newtheorem{conj}[thm]{Conjecture}
\newtheorem{qn}{Question}
\newtheorem{pro}{Proposition}[section]
\theoremstyle{definition}
\newtheorem{defn}{Definition}[section]
\newtheorem{remk}{Remark}[section]
\newcommand\bC{\mathbb{C}}
\newcommand\bR{\mathbb{R}}
\newcommand\bZ{\mathbb{Z}}
\def\square{\hfill${\vcenter{\vbox{\hrule height.4pt \hbox{\vrule width.4pt
height7pt \kern7pt \vrule width.4pt} \hrule height.4pt}}}$}

\newenvironment{pf}{{\it Proof.}\quad}{\square \vskip 12pt}

\title[Branched Spherical CR structures] {Branched Spherical CR structures on the complement of the figure eight knot.}

\author{Elisha Falbel}
\address{Institut de Math\'ematiques de Jussieu\\ 
Unit\'e Mixte de Recherche 7586 du CNRS \\
Universit\'e Pierre et Marie Curie \\
4, place Jussieu 75252 Paris Cedex 05, France
}
\author{Jieyan Wang}
\address{College of Mathematics and Econometrics,
 Hunan University, Changsha, 410082, People's Republic of China.} 
\thanks{The second author is  grateful to  Universit\'e Pierre et Marie Curie for the hospitality and support and to CSC and NSF(No. 11071059) for financial support. He also thanks Yueping Jiang, Baohua Xie and Wenyuan Yang for their encouragements.}

\begin{abstract}
We obtain a branched spherical CR structure on the complement of the figure eight knot whose holonomy representation
was given in \cite{fal08}.   There are essentially two boundary unipotent representations from the complement  of the figure eight knot into $\mathbf{PU}(2,1)$, we call them $\rho_1$ and $\rho_2$.  We make explicit some fundamental differences between these two representations.  For instance,  seeing the 
figure eight knot complement as a surface bundle over the circle,  the behaviour of of the
fundamental group of the fiber under the representation is a key difference between $\rho_1$ and $\rho_2$.

\end{abstract}

\maketitle

\section{Introduction}

The three dimensional sphere contained in $\bC^2$ inherits a Cauchy-Riemann structure as the boundary of
the complex two-ball.  Three dimensional manifolds locally modeled on the sphere then are called
spherical CR manifolds and have been studied since Cartan (\cite{C}).
Spherical CR structures  appear naturally as quotients of an open subset of the three dimensional
sphere by a subgroup of the CR automorphism group (denoted $\mathbf{PU}(2,1)$) (see \cite{J,Gol} and \cite{S} for a recent introduction).

The irreducible representations of the fundamental group  of the complement of the figure eight knot 
into $\mathbf{PU}(2,1)$ with unipotent boundary holonomy were obtained in \cite{fal08}.    To obtain such representations, one imposes the existence of a developing
map obtained from the 0-skeleton of an ideal triangulation.  Solution of a system of algebraic equations gives rise
to a set of representations of  $\Gamma= \pi_1(M)$, the fundamental group of the complement of the figure eight knot
with parabolic peripheral group.  

Up to pre-composition with automorphisms of $\Gamma$
there exists 2 irreducible representations into $\mathbf{PU}(2,1)$ with unipotent boundary holonomy
(see \cite{DF}). Following \cite{fal08} we call them $\rho_1$ and $\rho_2$.
In \cite{fal08} we showed that   $\rho_1$ could be obtained from a branched
spherical CR structure on the knot complement.  Moreover, this representation is not the holonomy of a complete structure
as the limit set is the full sphere $S^3$.

In this paper we analyze   $\rho_2$ and show that it
is also obtained as the holonomy of a branched structure in Theorem \ref{maintheorem} in section \ref{branched}.
The proof consists of extending the developing map obtained from the 0-skeleton to a developing map defined on simplices.  A complete (non-branched) spherical CR structure on the complement of the figure eight knot was obtained in \cite{DF}. Although the complete structure with unipotent boundary holonomy is unique (see \cite{DF}), it is not clear to us how to describe all branched structures.   The motivation to study branched CR structures is the hope that they would be 
easier to associate to a manifold once a representation is given.  As we have a general method to construct
representations of the fundamental group into  $\mathbf{PU}(2,1)$ we would like  an efficient method to 
obtain spherical CR structures with holonomy the given representation.  Constructing branched structures might be a step in this process.
Remark that a complete structure in the Whitehead link complement is described in  \cite{S} and, more recently, a whole family in \cite{PW}.

Another motivation for this paper is to stress a major difference between the two representations $\rho_1$ and $\rho_2$.  
Recall that the fundamental group of the figure eight knot complement contains a surface group (a punctured torus group) as a normal subgroup corresponding to the fundamental group of the fiber of the fibration of the complement over a circle. 
 In fact the kernel of the first representation  is contained in the surface group (and is not finitely generated)
 but the kernel of the second one is not.   This, in turn, implies that
the image of the surface group is of infinite index in the image of $\rho_1$  but of finite index in the image of  $\rho_2$. 
 Both images of the representations are contained in arithmetic lattices
as infinite index subgroups.  It turns out that the limit set of the image of  $\rho_1$  is 
the full $S^3$ but  the image of  $\rho_2$ has a proper limit set (see \cite{DF}).  These properties are given in sections \ref{rho1} and \ref{rho2}.   They might be general properties of representations of 3-manifold groups into $\mathbf{PU}(2,1)$.

We thank M. Deraux, A. Guilloux,  A. Reid,  P. Will and M.  Wolff for fruitful discussions.

\section{Complex hyperbolic space and its boundary.}
In this section, we introduce some basic materials about complex hyperbolic geometry. We refer Goldman's book \cite{Gol} for details.

\subsection{Complex hyperbolic space and its isometry group.}
Let $\mathbb{C}^{2,1}$ be the three dimensional complex vector space equipped with the Hermitian form
$$
\langle Z, W \rangle=Z_1\overline{W}_3+Z_2\overline{W}_2+Z_3\overline{W}_1.
$$
One has three subspaces:
$$
  V_{+}=  \{ Z\in \mathbb{C}^{2,1} : \langle Z, Z \rangle >0 \},
$$
$$
   V_{0}=  \{ Z\in \mathbb{C}^{2,1}- \{0\} : \langle Z, Z \rangle =0 \},
$$
$$
V_{-}=  \{ Z\in \mathbb{C}^{2,1} : \langle Z, Z \rangle < 0 \}.
$$
Let $P : \mathbb{C}^{2,1}- \{0\} \rightarrow \mathbb{C}P^2 $
be the canonical projection onto complex projection space.
Then  complex hyperbolic 2-space is defined as $\mathbf{H}^2_{\mathbb{C}}=P(V_{-})$ equipped with the Bergman metric.
The boundary of complex hyperbolic space is defined as $\partial \mathbf{H}^2_{\mathbb{C}}=P(V_{0}) $.

Let $\mathbf{U}(2,1)$ be the linear matrix group preserving the Hermitian form $\langle .,. \rangle$. The holomorphic isometry group $\mathbf{PU}(2,1)$ of $\mathbf{H}^2_{\mathbb{C}}$ is the projection of the unitary group $\mathbf{U}(2,1)$. The isometry group of $\mathbf{H}^2_{\mathbb{C}}$ is
$$
\widehat{\mathbf{PU}(2,1)}=\langle \mathbf{PU}(2,1), Z\mapsto \overline{Z} \rangle,
$$
where $Z\mapsto \overline{Z}$ is the complex conjugation.

The elements of $ \mathbf{PU}(2,1)$ can be classified into three
kinds of classes. Any element $g\in \mathbf{PU}(2,1)$ is called loxodromic if $g$ fixes exactly two points in $\partial \mathbf{H}^2$;
$g$ is called parabolic if it fixes exactly one point in $\partial \mathbf{H}^2$; otherwise, $g$ is called elliptic.

\subsection{Lattices.}

Let ${\mathcal O}_d$ be the ring of integers in the imaginary quadratic number
field ${\mathbb Q}(i\sqrt{d})$ where $d$ is a positive square-free integer.
If $d\equiv 1,2$ (mod 4) then ${\mathcal O}_d=\bZ[i\sqrt{d}]$ and if $d\equiv 3$ (mod 4) then ${\mathcal O}_d=\bZ[\frac{1+i\sqrt{d}}{2}]$.
The subgroup of ${\mathbf{PU}(2,1)}$ with entries in ${\mathcal O}_d$ is called the
{\sl Picard modular group} for ${\mathcal O}_d$ and is written
${\mathbf{PU}}(2,1;{\mathcal O}_d)$.  They are arithmetic lattices first considered by Picard.

\subsection{Heisenberg group and $\mathbb{C}-$circles}
The Heisenberg group $\mathfrak{N}$ is defined as the set $\mathbb{C}\times \mathbb{R}$ with group law
$$
(z,t)\cdot (z',t')=(z+z',t+t'+2\rm{Im}(z\overline{z}')).
$$
The boundary of complex hyperbolic space $\partial \mathbf{H}^2_{\mathbb{C}}$ can be identified with the
one point compactification $\overline{\mathfrak{N}}$ of $\mathfrak{N}$.

A point $p=(z,t)\in \mathfrak{N}$ and the point at infinity are lifted to the following points in $\mathbb{C}^{2,1}$:
$$
\hat{p}=\left[\begin{array}{c}
          (-|z|^2+it)/2 \\
          z \\
          1
        \end{array}\right]
\quad \rm{and} \quad
\hat{\infty}=\left[\begin{array}{c}
               1 \\
               0 \\
               0
             \end{array}\right].
$$

There are two kinds of totally geodesic submanifolds of real dimension 2 in $\mathbf{H}^2_{\mathbb{C}}$: complex geodesics and
totally real totally geodesic planes. Their boundaries in $\partial \mathbf{H}^2_{\mathbb{C}}$ are called $\mathbb{C}-$circles and
$\mathbb{R}-$circles.  Complex geodesics can be parametrized by their polar vectors, that is, points in $P(\mathbb{C}^{2,1})$
which are projections of vectors orthogonal to the lifted complex geodesic.

\begin{pro}
In the Heisenberg model,  $\mathbb{C}-$circles are either vertical lines or ellipses, whose projection on the z-plane are circles.
\end{pro}

For a given pair of distinct points in $\partial \mathbf{H}^2_{\mathbb{C}}$, there is a unique $\mathbb{C}-$circle passing through them.
Finite $\mathbb{C}-$circles are determined by a centre and a radius. For example, the finite $\mathbb{C}-$circle with centre $(z_0,t_0)$ and
radius $R>0$ has polar vector
$$
\left[
  \begin{array}{c}
    (R^2-|z_0|^2+it_0)/2 \\
    z_0 \\
    1 \\
  \end{array}
\right]
$$
and in which any point $(z,t)$ satisfies the equations
$$
\left\{
  \begin{array}{ll}
    |z-z_0|=R \\
    t=t_0+2\rm{Im}(\overline{z}z_0)
  \end{array}
\right.
$$
\subsection{CR structures}\label{sec:cr}
CR structures appear naturally as boundaries of complex manifolds.  The local geometry
of these structures was studied by E. Cartan \cite{C} who defined, in dimension three, a curvature analogous to
curvatures of a Riemannian structure.  When that curvature is zero,  
Cartan called them spherical CR structures and developed their basic properties.  A much later  study by
Burns and Shnider \cite{BS} contains the modern setting for these structures. 

\begin{defn}
A spherical CR-structure on a 3-manifold is a geometric structure
modeled on the homogeneous space $S^3$ with the above
${\mathbf{PU}(2,1)}$ action.
\end{defn}

\begin{defn}
We say a spherical CR-structure on a 3-manifold is complete if it is
equivalent to a quotient of the domain of regularity in $S^3$ by a
discrete subgroup of ${\mathbf{PU}(2,1)}$.
\end{defn}

Here, equivalence between CR structures is defined, as usual, by
diffeomorphisms preserving the structure.  The diffeomorphism group
of a manifold therefore acts trivially on its CR structures.  
Observe that taking the complex conjugate of local charts of a CR structure (maps of open sets 
into $S^3$) gives another CR structure which might not be equivalent to the original one.   A weaker definition
of spherical CR structures as geometric structures modeled on 
 the full isometry group $\widehat{{\mathbf{PU}(2,1)}}$ 
with its action on $S^3$ is sometimes preferable.  Indeed, in that case,  complex conjugation of
local charts will induce an equivalent spherical CR structure.

A CR structure, in particular, has an orientation which is compatible
with the orientation induced by its contact structure.  Observe that even $\widehat{{\mathbf{PU}(2,1)}}$ preserves orientation
so a spherical CR structure in the weaker sense is also oriented.   Both
orientations of $S^3$ are obtained via equivalent CR structures
because there exists an orientation reversing diffeomorphism of $S^3$.
More generally, manifolds which have orientation reversing maps either
have equivalent CR structures opposite orientations or none.  On the
other hand, It is not clear if a manifold having a CR structure will
have another one giving its opposite orientation.%

As all geometric structures, a spherical structure on a manifold $M$
induces a developing map defined on its universal cover $\tilde M$
$$                                                                                                                            
d : \tilde M \rightarrow S^3                                                                                                 
$$
and a holonomy representation
$$                                                                                                                            
\rho : \pi_1(M)\rightarrow{\mathbf{PU}(2,1)}.                                                                                  
$$

Observe again that pre-composition with a diffeomorphism will induce an equivalent structure
with a holonomy representation which is obtained from the old one by pre-composition with
an automorphism of the fundamental group. 
Also observe that the holonomy representation is not discrete in general and the developping map
might be surjective.  


\subsection{Branched structures}

Given a representation is not clear that it is defined as the holonomy representation of a spherical 
CR structure.  In that sense it is useful to introduce a weaker definition of branched structure in the hope that
representations might be understood in a geometric way.

A branched spherical CR structure is a CR structure except along some curves where the structure is locally modeled on the $t$-axis
inside $\bR^3= \{ \ (z,t)\     | \ z\in \bC, t\in \bR \ \}$ together with the ramified map into the Heisenberg group given by
$$
(z,t)\rightarrow (z^n,t),
$$
where $n$ is the branching order.  The CR structure around the curve is given by the pullback of the CR structure around 
the Heisenberg $t$-axis.

\section{The figure eight knot.}

We use the same notations as that in the paper \cite{fal08} and recall briefly the three irreducible representations obtained there.

The figure eight knot complement $M$ has a fundamental group $\Gamma= \pi_{1}(M)$ which can be presented as
$$
\Gamma=\langle \ g_1, g_3\  | \ [g_3,g_1^{-1}]g_3= g_1[g_3,g_1^{-1}]\ \rangle.
$$
It is useful to introduce another generator
$$
g_2 = [g_3,g_1^{-1}],
$$
that is
$$
g_1 = g_2g_3g_2^{-1}.
$$
The figure eight knot complement is fibered over the circle with fiber a punctured torus.  The fibration is
encoded in the following sequence.
$$
1\rightarrow F_2\rightarrow \Gamma\rightarrow \bZ\rightarrow 0
$$
Here, $F_2$ is the free group of rank 2 with generators
$$F_2=\langle \ a=g_2,b=[g_2,g_3^{-1}] \ \rangle.
$$
We can then present
$$
\Gamma = \langle \ a, b, t\  | \ tat^{-1}= aba , tbt^{-1}= ab \ \rangle.
$$
where $t=g_3$ is seen to act as a pseudo-Anosov element of the mapping class group of $F_2$.

We consider in this paper the following representations into $\mathbf{SU}(2,1)$ obtained in \cite{fal08}:
\begin{enumerate}
\item $$
 \rho_1(g_1)=\begin{pmatrix}
           1&  1 &     -\frac{1}{2}- \frac{\sqrt{3}}{2}i\\
 0 & 1 & -1  \\
0 & 0 & 1\\
           \end{pmatrix},
\  \
           \rho_1( g_3)=\begin{pmatrix}
           1&  0 &   0\\
 1 & 1 & 0  \\
  -\frac{1}{2}- \frac{\sqrt{3}}{2}i & -1 & 1\\
           \end{pmatrix}.
           $$

\item
 $$
 \rho_2(g_1)=\begin{pmatrix}
           1&  1 &   -\frac{1}{2}- \frac{\sqrt{7}}{2}i\\
 0 & 1 & -1  \\
0 & 0 & 1\\
           \end{pmatrix},
 \ \
            \rho_2(g_3)=\begin{pmatrix}
           1&  0 &   0\\
 -1 & 1 & 0  \\
-\frac{1}{2}+ \frac{\sqrt{7}}{2}i& 1 & 1\\
           \end{pmatrix}.
           $$

\item
$$
\rho_3( g_1)=\begin{pmatrix}
           1&  1 &   -1/2\\
 0 & 1 & -1  \\
0 & 0 & 1\\
           \end{pmatrix},
\ \             \rho_3(g_3)=\begin{pmatrix}
           1&  0 &   0\\
\frac{5}{4}- \frac{\sqrt{7}}{4}i  & 1 & 0  \\
-1& -\frac{5}{4}- \frac{\sqrt{7}}{4}i & 1\\
           \end{pmatrix}.
           $$
\end{enumerate}

The representation $\rho_3$ is obtained by pre-composition of $\rho_2$ with the automorphisms of the fundamental group associated
to a reversing orientation diffeomorphism.  So we will concentrate in the first two representations during the rest of this paper.

\section{Representations}

As every complement of a tame knot, the complement of the figure eight knot has fundamental group $\Gamma$ fitting in the
exact sequence
$$
1\rightarrow [\Gamma,\Gamma]\rightarrow \Gamma\rightarrow \bZ\rightarrow 0.
$$
In the case of the complement of the figure eight knot we have
$$
1\rightarrow F_2\rightarrow \Gamma\rightarrow \bZ\rightarrow 0
$$
where $F_2$ is the free group of rank two.  We will be interested in the general case when
$$
1\rightarrow F \rightarrow \Gamma\rightarrow \bZ\rightarrow 0
$$
is an exact sequence.
Suppose
$$
\rho : \Gamma\rightarrow G
$$
is a representation with $K=Ker (\rho)$.

\begin{lem}  \label{lemma:K}The following diagram is commutative:
$$
\begin{CD}
 @. 1 @. 1 @. 1 @.\\
@. @V VV @V VV @ V VV@. \\
1 @> >> K\cap  F @> >> K @> p >> p(K) @> >>  0\\
@. @V VV @V VV @ V VV@. \\
1 @> >>  F @> >> \Gamma @> p >> \bZ @> >>  0\\
@. @V\rho VV @V\rho VV @ V\bar \rho VV@. \\
1 @> >>\rho( F) @> >> \rho(\Gamma) @> \bar p >> \rho(\Gamma)/\rho(F)@> >>  0\\
@. @V VV @V VV @ V VV@. \\
 @. 1 @. 1 @. 1 @.
\end{CD}
$$
Where $\bar p$ is the quotient map and $\bar \rho$ is defined so that the diagram be commutative.
\end{lem}

\begin{pf} The only verification we have to make is that $Ker(\bar \rho)$ is the image of $p(K)$.  Suppose
$x=p(ft^n)\in Ker(\bar \rho)$  with $t\in p^{-1}(1)$, $n\in \bZ$ and $f\in F$ satisfying
$$
\bar p \rho( ft^n) =Id.
$$
Then $\rho( ft^n)=\rho( f')$ with $f'\in  F$.  Therefore ${f'}^{-1}ft^n\in  K$ and then
$x=p(ft^n)=p({f'}^{-1}ft^n)\in p(K)$.
\end{pf}

We conclude that the inclusion $\rho( F)\subset  \rho(\Gamma)$ is of finite index if and only if
$K$ contains an element $f t^n$,$n\neq 0$, where $f\in F$ and $p(t)=1$.  The index is precisely the least absolute value of an integer
satisfying the condition.

\begin{cor} \label{infinite}
$\rho( F)\lhd  \rho(\Gamma)$ is of infinite index if and only if $K\subset F$.
\end{cor}

\subsection{}\label{rho1}
\begin{center}
\bf{The representation} $\rho_1$.
\end{center}

We consider the first representation. Let $\omega_3=-\frac{1}{2}+i\frac{\sqrt{3}}{2}$.
The ring of integers of the field $Q(i \sqrt{3})$ is ${\mathcal O}_3=\mathbb{Z}[\omega_3]$.
 The representation is discrete, since the generators
$G_1,G_2,G_3$ are contained in the arithmetic lattice $P_3=PU(2,1;{\mathcal O}_3)$.

We use the presentation of $P_3$ obtained in \cite{FP}:
$$
P_3=\langle P, Q, I\ |\ I^2=(QP^{-1})^6=PQ^{-1}IQ P^{-1}I=P^3Q^{-2}=(IP)^3\ \rangle.
$$

Recall from \cite{fal08}
that $\rho_1 (\Gamma)$ is generated by
$$
G_1= [P,Q]\ ,
$$
$$
G_2= [I,[Q,P]]\ ,
$$
$$
G_3= A[P,Q]A^{-1}\ ,
$$
with $A=P^{-2}IP^2$.

A usefull tool  in the following computations is the normalizer $N=N(\rho_{1} (\Gamma))\subset P_3$, the least normal subgroup of $P_3$ containing $\rho_1 (\Gamma)$.

\begin{lem} [\cite{fal08}]
$P_3/N$ is isomorphic to $\bZ/6\bZ$.
\end{lem}

Computing that $P_3/[P_3,P_3]$ is of order 6  and observing that $[P_3,P_3]\subset N$ we remark that $N=[P_3,P_3]$.
By computing the quotient of $P_3$ by the normalizer of $\langle G_1 , I, [Q,P^{-1}] \rangle $ we obtain the following
\begin{lem} [\cite{fal08}]
$N=\langle G_1 , I, [Q,P^{-1}] \rangle .$
\end{lem}

\begin{lem}
$P_3/[N,N]$ is isomorphic to the euclidean triangle group of type $(2,3,6)$.
\end{lem}
\begin{pf}
Using the presentation of $P_3$ and the lemma above we obtain for the presentation of the quotient
$$
P_3/[N,N]=\langle P, Q\ |\ (QP^{-1})^6=P^3=Q^{2}\ \rangle.
$$
\end{pf}

\begin{lem}
$N/[N,N]$ is isomorphic to $\bZ\oplus \bZ$.
\end{lem}
\begin{pf}
From the isomorphism theorem
$$
N/[N,N]=\frac{P_3/[N,N]}{P_3/N}
$$
and the two previous lemmas we obtain the result.
 \end{pf}

We have (cf. \cite{fal08})
$$
[\rho_1(\Gamma),\rho_1(\Gamma)]\lhd \rho_1(\Gamma)
\lhd \langle G_1, I \rangle
\lhd N
\lhd
PU(2,1;{\mathcal O}_3)
$$
with the last inclusion of order 6  and the inclusion $\rho_1(\Gamma)
\lhd \langle G_1, I \rangle$ of index at most two.

Observe now that the inclusion  $\langle G_1, I \rangle
\lhd N $ has abelian quotient and therefore
$[N,N]\subset \langle G_1, I \rangle$ so we obtain
$$
[\rho_1(\Gamma),\rho_1(\Gamma)]\subset [N,N]
\lhd \langle G_1, I \rangle
\lhd  N
\lhd
PU(2,1;{\mathcal O}_3)
$$
The following Proposition was obtained after discussions with A. Reid.  The proof given here is a simplification 
of his argument which involved a gap computation (\cite{R}).
\begin{pro}\label{proposition:finite}
The inclusions
$$
[\rho_1(\Gamma),\rho_1(\Gamma)]\lhd \rho_1(\Gamma)
\subset
PU(2,1;{\mathcal O}_3)
$$
are of infinite index.
\end{pro}
\begin{pf}
Observe first that $[\rho_1(\Gamma),\rho_1(\Gamma)]\lhd  \langle G_1, I \rangle$ and $[N,N]\lhd  \langle G_1, I \rangle$ are two normal inclusions
and therefore
$$
 \langle G_1, I \rangle/[N,N]\rightarrow  \langle G_1, I \rangle/[\rho_1(\Gamma),\rho_1(\Gamma)]
 $$
 is a monomorphism.  On the other hand, the quotient $\langle G_1, I \rangle/[\rho_1(\Gamma),\rho_1(\Gamma)]$ is finite or contains $\bZ$ as a subgroup of index at most two.

 Suppose now that $\rho_1(\Gamma)
\subset
PU(2,1;\mathbb{Z}[\omega])$ is of finite index.  Then $\langle G_1, I \rangle\lhd N$ should be of finite index and therefore, as $N/[N,N]=\bZ\oplus \bZ$, $\langle G_1, I \rangle/[N,N]=\bZ\oplus \bZ$.
This contradicts the monomorphism above.

Suppose next that $[\rho_1(\Gamma),\rho_1(\Gamma)]\lhd \rho_1(\Gamma)$ is of finite index.  Then the inclusion $[\rho_1(\Gamma),\rho_1(\Gamma)]\lhd \langle G_1, I \rangle$ would be of finite index.
This, in turn,  implies that $[N,N]\lhd \langle G_1, I \rangle$ is of finite index.  Now this contradicts the monomorphism
$$
 N/[N,N]\rightarrow  N/\langle G_1, I \rangle
 $$
as $N/\langle G_1, I \rangle$ is abelian of rank at most one.
\end{pf}

From Lemma \ref{lemma:K} and Proposition \ref{proposition:finite}
we obtain the following
\begin{cor}
$$
Ker(\rho_1) \lhd [\Gamma,\Gamma ].
$$
\end{cor}
  We conclude with the following property of the kernel:

\begin{pro}
$Ker(\rho_1)$ is not of finite type.
\end{pro}
\begin{pf}
Observe that  $Ker(\rho_1)$ is  clearly preserved under the pseudo-anosov element of the mapping class group denoted by $t$.
The result then follows from Lemme 6.2.5 in \cite{O}.
\end{pf}

\subsection{}\label{rho2}
\begin{center}
{\bf{The representation} $\rho_2$.}
\end{center}


The second representation (see the subsection 6.5.1 in \cite{fal08}) is given by $\Gamma_2=\rho_2(\pi_1(M))$, with $\Gamma_2=\langle \rho_2(g_1),\rho_2(g_2), \rho_2(g_3) \rangle$, where
$$
G_1:=\rho_2(g_1)=\left(
  \begin{array}{ccc}
    1 & 1 & -\frac{1}{2}-i\frac{\sqrt{7}}{2} \\
    0 & 1 & -1 \\
    0 & 0 & 1 \\
  \end{array}
\right)
$$
$$
G_2:=\rho_2(g_2)=\left(
  \begin{array}{ccc}
    2 & \frac{3}{2}-i\frac{\sqrt{7}}{2} & -1 \\
    -\frac{3}{2}-i\frac{\sqrt{7}}{2} & -1 & 0 \\
    -1 & 0 & 0 \\
  \end{array}
\right)
$$
$$
G_3:=\rho_2(g_3)=\left(
  \begin{array}{ccc}
    1 & 0 & 0 \\
    -1 & 1 & 0 \\
    -\frac{1}{2}+i\frac{\sqrt{7}}{2} & 1 & 1 \\
  \end{array}
\right).
$$
Moreover, $G_2=[G_3,G_1^{-1}]$ is a regular elliptic element of order four, and $G_1,G_3$ are pure parabolic elements.
\begin{remk}
We can also see that the element $G_3G_1^{-1}$ is loxodromic. The fixed points of $G_1, G_3$ are respectively
$p_1=\infty$ and $p_2=(0,0)$.
\end{remk}

Let $\omega_7=\frac{1}{2}+i\frac{\sqrt{7}}{2}$.
The ring of integers of the field $Q(i \sqrt{7})$ is ${\mathcal O}_7=\mathbb{Z}[\omega_7]$.  We observe then that
 the representation is discrete, since the generators
$G_1,G_2,G_3$ are contained in the arithmetic lattice $PU(2,1; {\mathcal O}_7)$.

\begin{thm}[see  Proposition 3.3 and Theorem 4.4 in \cite{zh12}]
The group $PU(2,1;{\mathcal O}_7)$ is generated by the elements
$$
I=\left(
    \begin{array}{ccc}
      0 & 0 & 1 \\
      0 & -1 & 0 \\
      1 & 0 & 0 \\
    \end{array}
  \right),
R_1=\left(
      \begin{array}{ccc}
        1 & 0 & 0 \\
        0 & -1 & 0 \\
        0 & 0 & 1 \\
      \end{array}
    \right),
R_2=\left(
      \begin{array}{ccc}
        1 & 1 & -\bar{\omega_7} \\
        0 & -1 & 1 \\
        0 & 0 & 1 \\
      \end{array}
    \right),
$$
$$
R_3=\left(
      \begin{array}{ccc}
        1 & \bar{\omega_7} & -1 \\
        0 & -1 & \omega_7 \\
        0 & 0 & 1 \\
      \end{array}
    \right),
T=\left(
    \begin{array}{ccc}
      1 & 0 & i\sqrt{7} \\
      0 & 1 & 0 \\
      0 & 0 & 1 \\
    \end{array}
  \right).
$$
Moreover, the stabilizer subgroup of  infinity  has the presentation
$$
\langle R_1,R_2,R_3,T | R_1^2=R_3^2=[T,R_1]=[T,R_3]=T R_2^{-2}=(R_1R_3R_2)^2=Id  \rangle.
$$
\end{thm}

We may express the generators of  $\rho_2(\Gamma)$ in terms of the generators of $PU(2,1; {\mathcal O}_7)$:
\begin{pro}
$$
G_1=R_1R_2T^{-1}=R_1R_2^{-1},
$$
$$
-G_2=R_2R_1R_3I,
$$
$$
G_3=IR_2IR_1=IR_2R_1I=IG_1^{-1}I.
$$
\end{pro}

We also observe that
$$
\rho_2(\Gamma)
=\langle G_1, G_2, G_3 \rangle
=\langle G_1, G_3 \rangle
\lhd \langle G_1, I \rangle
$$
where $\langle G_1, G_3 \rangle \lhd \langle G_1, I \rangle$ is a subgroup of index at most two since
$G_3=IG_1^{-1}I$.

We also have
\begin{lem}
$$
\langle G_1, I, T \rangle \lhd \langle G_1, I, T, R_1 \rangle
=PU(2,1;{\mathcal O}_7)
$$
\end{lem}
\begin{pf}
$ \langle G_1, I, T \rangle \lhd \langle G_1, I, T, R_1 \rangle$ is a normal subgroup since
$R_1G_1R_1=T^{-1}G_1^{-1}$, $R_1IR_1=I$ and $R_1TR_1=T$.  The normal inclusion is of index at most two.
\end{pf}

The inclusion
 $$\langle G_1, I \rangle
\subset \langle G_1, I, T \rangle$$
 can be neither normal nor finite if one proves that the limit set is not $S^3$.

A simple computation shows that
\begin{lem}$$
 \rho_2(t^3)=[\rho_2(a^{-1}),\rho_2(b^{-1})]
 $$
\end{lem}

From the lemma above and Lemma \ref{infinite} we obtain
\begin{cor}
$$
[\rho_2(\Gamma),\rho_2(\Gamma)]\lhd \rho_2(\Gamma)
$$
is of index at most three.
\end{cor}

\section{Tetrahedra}

\subsection{Edges}
\begin{figure}[h!!]
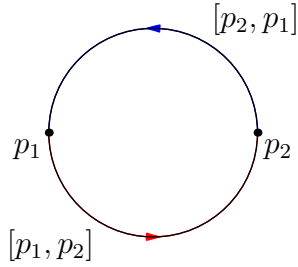
 \label{circle}
\begin{center}
\begin{asy}
  // 
  // 
  size(4cm,0);
  //defaultpen(1);
 // usepackage("amssymb");
  //import geometry;
// On dŽfinit un rŽel et trois points.
real h=5;
pair O=(0,0),A=(-h,0),B=(h,0);
// On trace le demi-cercle de diamtre [AB]
// en tournant dans le sens direct de B vers A.
draw(arc(O,B,A),blue,MidArrow(5bp));draw(arc(O,A,B),red,MidArrow(5bp));
// Et on ajoute les points.
dot("$p_1$",A,SW);
dot("$p_2$",B,SE);

draw(Label("$[p_2,p_1]$",Relative(1/3),align=NE),arc(O,B,A));
draw(Label("$[p_1,p_2]$",Relative(1/3),align=SW),arc(O,A,B));
\end{asy}
\caption{A $\bC$-circle between two points in $S^3$ is oriented and defines two oriented segments.}
\end{center}
\end{figure}

Given two points $p_1$ and $p_2$ in $S^3$, there exists a unique $\bC$-circle between them.  As the boundary of a complex disc
has a positive orientation, the $\bC$-circle inherits that orientation and defines therefore two distinct arcs $[p_1,p_2]$ and $[p_2,p_1]$ (see Figure \ref{circle}).

The complex disc is contained in complex hyperbolic space and does not intersect $S^3$ except in the boundary.  If one wants to obtain a disc in $S^3$ whose boundary is a $\bC$-circle, a usefull construction is
obtained by using a family of $\bC$-circles which foliates the disc. It could have a singularity at one point in the interior or at the boundary.  Remark though that this construction
is not canonical.

\subsection{Triangles}

Given three points  $p_1,p_2,p_3\in S^3$ we might construct six different triangles (1-skeletons) corresponding to an orientation choice of the edges between the points.  Observe that if the three points are contained in the same $\bC$-circle some edges contain two vertices.  If we suppose, on the other hand, that the three points are not in the same $\bC$-circle (we refer to it as a generic configuration) then the three edges (any choice) intersect only at the vertices.

In order to obtain a surface whose boundary is the triangle we might fix one point $p\in S^3$ and  consider the segments (of $\bC$-circles) joining that point and the edges as in a barycentric construction.   If the orientation of the edges permits, one can degenerate this construction making $p$ approach one of the vertices (it is clear that this is not  possible only if the edges define an orientation of the triangle).  There are choices to be made in that construction and each choice corresponds to a different triangle.

It is easier to analyse first the case where the triangle is degenerate:  If the vertices are generic the triangle defined above is embedded.  In the case of an oriented triangle, we could add first
an edge complement of one of the edges (making a full $\bC$-circle) and then consider a degenerate triangle with that edge union a disc whose boundary is the full $\bC$-circle now foliated by
$\bC$-circles with singular point one of the vertices.

It is not clear from the definition that the surface  defined is embedded.  Each case needs a verification.

\begin{figure}[h!!!]
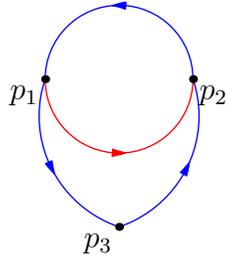
 \label{triangle}
\begin{center}
\begin{asy}
  // 
  // 
  size(3cm,0);
  //defaultpen(1);
 // usepackage("amssymb");
  //import geometry;
 real h=5;
path e2=(-h,0)..(-h,-5)..(0,-2*h);
path e1=(0,-2*h)..(h,-5)..(h,0);
draw(e2,blue,MidArrow(5bp));draw(e1,blue,MidArrow(5bp));

// On dŽfinit un rŽel et quatre points.
pair O=(0,0),A=(-h,0),B=(h,0),D=(0,-2*h);
// On trace le demi-cercle de diamtre [AB]
// en tournant dans le sens direct de B vers A.
draw(arc(O,B,A),blue,MidArrow(5bp));draw(arc(O,A,B),red,MidArrow(5bp));

// Et on ajoute les points.
dot("$p_1$",A,SW);
dot("$p_2$",B,SE);
dot("$p_3$",D,SW);
//draw(Label("$[p_2,p_1]$",Relative(1/3),align=NE),arc(O,B,A));
//draw(Label("$[p_1,p_2]$",Relative(1/3),align=SW),arc(O,A,B));
\end{asy}
\caption{Triangles.  We show four possible edges given a configuration of three points.}
\end{center}
\end{figure}

\begin{figure}[h!!!]
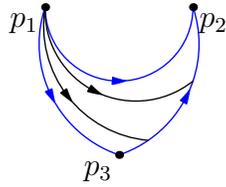
 \label{2-skeleton}
\begin{center}
\begin{asy}
  // 
  // 
  size(3cm,0);
  //defaultpen(1);
 // usepackage("amssymb");
  //import geometry;
 real h=5;
path e2=(-h,0)..(-h,-5)..(0,-2*h);
path e1=(0,-2*h)..(h,-5)..(h,0);
draw(e1,blue,MidArrow(5bp));draw(e2,blue,MidArrow(5bp));

// On dŽfinit un rŽel et quatre points.
pair O=(0,0),A=(-h,0),B=(h,0),D=(0,-2*h);
// On trace le demi-cercle de diamtre [AB]
// en tournant dans le sens direct de B vers A.
//draw(arc(O,B,A),blue,MidArrow(5bp));
draw(arc(O,A,B),blue,MidArrow(5bp));

// Et on ajoute les points.
dot("$p_1$",A,SW);
dot("$p_2$",B,SE);
dot("$p_3$",D,SW);

// Le 2-squelete
 pair u1=(-2,-h-.5);
pair u2=(-3.5,-h-1);

path uu=A..u1..point(e1,1);
path uu2=A..u2..point(e1,.3);
draw(uu,MidArrow(5bp));draw(uu2,MidArrow(5bp));

\end{asy}
\caption{Triangles.  A 2-skeleton foliated by arcs in $\bC$-circles.}
\end{center}
\end{figure}

\begin{figure}[h!!!]
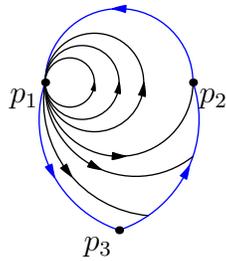
 \label{2triangle}
\begin{center}
\begin{asy}
  // 
  // 
  size(3cm,0);
  //defaultpen(1);
 // usepackage("amssymb");
  //import geometry;
 real h=5;
path e2=(-h,0)..(-h,-5)..(0,-2*h);
path e1=(0,-2*h)..(h,-5)..(h,0);
draw(e2,blue,MidArrow(5bp));draw(e1,blue,MidArrow(5bp));

// On dŽfinit un rŽel et quatre points.
pair O=(0,0),A=(-h,0),B=(h,0),D=(0,-2*h);
// On trace le demi-cercle de diamtre [AB]
// en tournant dans le sens direct de B vers A.
draw(arc(O,B,A),blue,MidArrow(5bp));draw(arc(O,A,B),MidArrow(5bp));

// Et on ajoute les points.
dot("$p_1$",A,SW);
dot("$p_2$",B,SE);
dot("$p_3$",D,SW);

// Le 2-squelete
 pair u1=(-2,-h-.5);
pair u2=(-3.5,-h-1);

path uu=A..u1..point(e1,1);
path uu2=A..u2..point(e1,.3);
draw(uu,MidArrow(5bp));draw(uu2,MidArrow(5bp));

draw(circle((-h+2*h/3,0),2*h/3),Arrow(5bp));draw(circle((-h+h/2,0),h/2),Arrow(4bp));draw(circle((-h+h/3,0),h/3),Arrow(3bp));

\end{asy}
\caption{Triangles.  A 2-skeleton foliated by arcs in $\bC$-circles. The case where the degenerate barycentric construction does not apply.}
\end{center}
\end{figure}

\subsection{CR Tetrahedra}

Once faces whose border are triangles are defined one can define a 3-simplex based on a configuration of four points by choosing faces
to each of the four configurations of three points.  The problem is that the choices have to be compatible and faces, otherwise well defined, could
intersect between one another.

We will make arguments
using sometimes flat discs adjoined to edges keeping in mind that we could, in fact, by a slight deformation deal with 3-simplices.

\begin{defn}
A tetrahedron is called a generalized tetrahedron if it has a disc adjoined to an edge.  That is is a simplex union a disc whose intersection with the simplex is an edge contained in the boundary of the disc.
\end{defn}

We could deform then these faces thickening the disc to obtain a topological 3-simplex.
Faces of tetrahedra are not canonical
and we will make use of this flexibility.

\section{Branched CR structures associated to  representations.}\label{branched}

\begin{figure}[h!!]
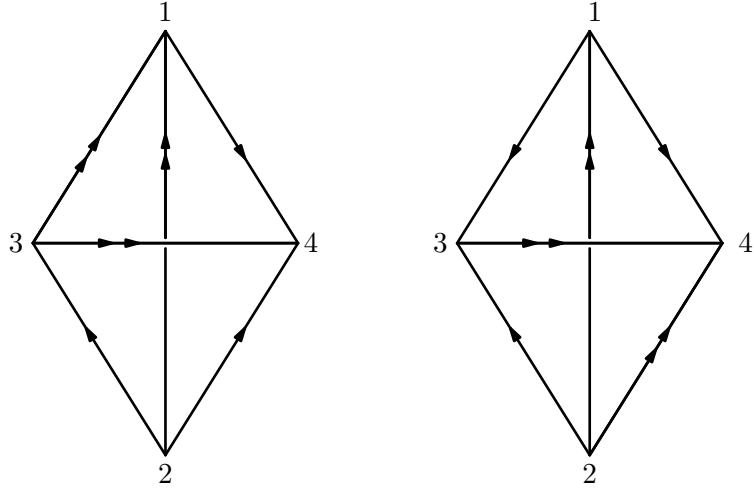
\label{eight}
\begin{center}
\begin{asy}
  // 
  // 
  size(10cm);
  defaultpen(1);
  usepackage("amssymb");
  import geometry;

  // 
  point o = (8,0);
  point oo = (24,0);
  pair ph = (10,0);
  pair pv = (5,8);

  // 
  draw(o -- o+pv,Arrow(5bp,position=.4), Arrow(5bp,position=.5));
  draw(o+pv -- o+ph, Arrow(5bp,position=.6));
  draw(o+ph-pv-- o, Arrow(5bp,position=.6));
  draw(o+ph-pv-- o+ph, Arrow(5bp,position=.6));
  draw(o+ ph/2+(0,.2)--o+pv, Arrow(5bp,position=.4), Arrow(5bp,position=.5));
  draw(o+ph-pv -- o+ph/2+ (0,-.2));
  draw(o -- o+ph,Arrow(5bp,position=.3), Arrow(5bp,position=.4) );

  label("{\small  $3$}", o, 1*W);

  label("{\small $4$}", o+ph, .6*E);

  label("{\small $1$}", o+pv, 1*N);
  label("{\small $2$}", o+ph -pv, 1*S);

  // 
  draw( oo+pv -- oo,Arrow(5bp,position=.6));
  draw( oo+pv -- oo+ph,Arrow(5bp,position=.6));
  draw( oo+ph-pv -- oo,Arrow(5bp,position=.6));
  draw(oo+ ph/2+(0,.2)--oo +pv,Arrow(5bp,position=.4), Arrow(5bp,position=.5));
  draw(oo+ph-pv -- oo+ ph/2-(0,.2));
  draw(oo+ph-pv -- oo+ ph,Arrow(5bp,position=.5),Arrow(5bp,position=.6));
  draw(oo -- oo+ph,Arrow(5bp,position=.3), Arrow(5bp,position=.4));

 label("{ \small $3$}", oo, 1*W);

   label("{ \small $4$}", oo+ph, .6*E);

  label("{ \small $1$}", oo+pv, 1*N);

  label("{\small $2$}", oo+ph -pv, 1*S);

\end{asy}
\caption{The figure eight represented by two tetrahedra with face pairings defined by the set of arrows.}
\end{center}
\end{figure}

The representations in \cite{fal08} are obtained by imposing that the 0-skeleton of an ideal triangulation defines a developing map.
The triangulation of the figure eight knot complement is shown in Figure \ref{eight}.  The 0-skeleton can be realized
as points in $S^3$ and using the side pairing maps we can define a developing map on the 0-skeleton of the universal covering.

In order to obtain a spherical CR structure we have to define the 1-skeleton, the 2-skeleton, then obtain 3-simplices and show that the developing
map defined on the 0-skeleton extends to the 3-simplices.

\begin{figure}[h!!]
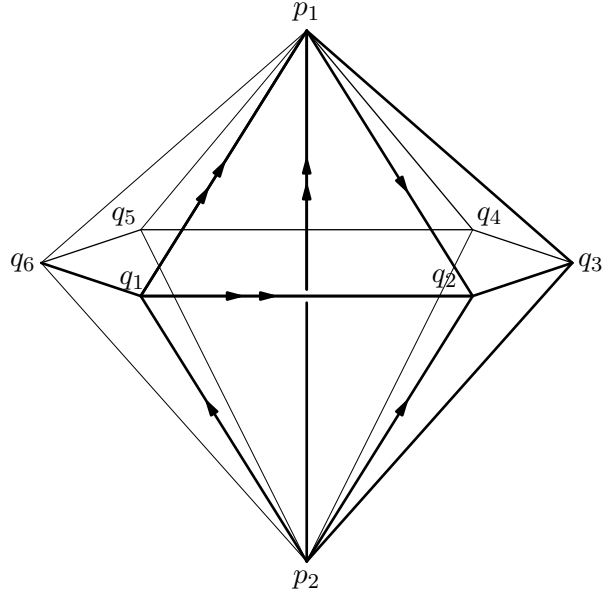
\label{alongedge}
\begin{center}
\begin{asy}
  // 
  // 
  size(8cm);
  defaultpen(1);
  usepackage("amssymb");
  import geometry;

  // 
  point o = (8,0);
  pair ph = (10,0);
  pair pv = (5,8);

  pair pq3 = (3,1); //
  pair q4 = (-3,1); //
   pair q5 = (-10,0); //
    pair q6= (-3,-1); //

  // 
  draw(o -- o+pv,Arrow(5bp,position=.4), Arrow(5bp,position=.5));
  draw(o+pv -- o+ph, Arrow(5bp,position=.6));
  draw(o+ph-pv-- o, Arrow(5bp,position=.6));
  draw(o+ph-pv-- o+ph, Arrow(5bp,position=.6));
  draw(o+ ph/2+(0,.2)--o+pv, Arrow(5bp,position=.4), Arrow(5bp,position=.5));
  draw(o+ph-pv -- o+ph/2+ (0,-.2));
  draw(o -- o+ph,Arrow(5bp,position=.3), Arrow(5bp,position=.4) );

   // 
  draw(o+pv -- o+ph+pq3);
  draw(o+ph-pv-- o+ph+pq3);
 draw(o+ph-- o+ph+pq3);

// 
 draw(o+ph+pq3-- o+ph+pq3 +q4,linewidth(.5bp));
draw(o+ph+pq3+q4-- o+ph+pq3 +q4+q5,linewidth(.5bp));
draw(o+ph+pq3+q4 +q5-- o+ph+pq3 +q4+q5+q6,linewidth(.5bp));
draw(o+ph+pq3+q4 +q5+q6-- o);

 draw(o+pv -- o+ph+pq3 +q4,linewidth(.5bp));
draw(o+pv -- o+ph+pq3 +q4+q5,linewidth(.2bp));
draw(o+pv -- o+ph+pq3 +q4+q5+q6,linewidth(.2bp));

 draw(o+ph-pv -- o+ph+pq3 +q4,linewidth(.2bp));
draw(o+ph-pv -- o+ph+pq3 +q4+q5,linewidth(.2bp));
draw(o+ph-pv -- o+ph+pq3 +q4+q5+q6,linewidth(.2bp));

  // 
  //
   label("{\small  $q_1$}", o, .7*dir(115));

  //
  label("{\small $q_2$}", o+ph, 2*dir(145));
   label("{\small $q_3$}", o+ph+pq3, .6*E);
    label("{\small $q_4$}", o+ph+pq3 +q4, .6*dir(45));
   label("{\small $q_5$}", o+ph+pq3+q4+q5, .6*NW);
    label("{\small $q_6$}", o+ph+pq3+q4+q5+q6, .6*W);

   //
   label("{\small $p_1$}", o+pv, 1*N);

   //
  label("{\small $p_2$}", o+ph -pv, 1*S);

\end{asy}
\caption{A schematic view of the six tetrahedra glued around the edge $[p_1,p_2]$.}
\end{center}
\end{figure}

Once we obtain two 3-simplices in $S^3$ which have well defined side pairings we might have some branching along  the edges of the
simplices.  In fact we will prove that around one of the edges the simplices are put together as in Figure \ref{alongedge} but along the other
edge we show that the 6 tetrahedra turn around the edge three times.

\begin{figure}[h!!]
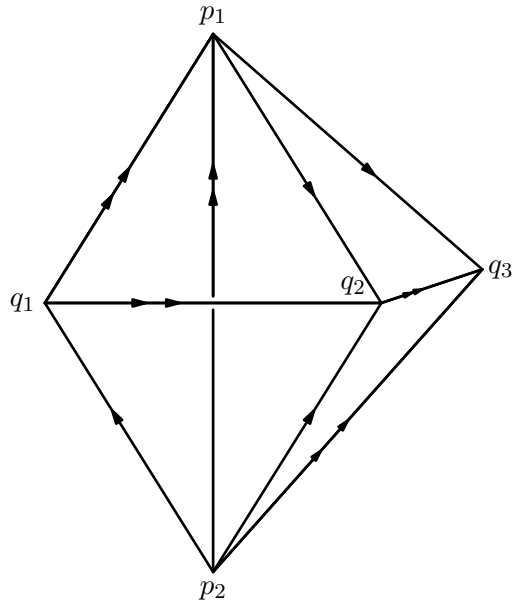

\begin{center}
\begin{asy}
  // 
  // 
  size(8cm);
  defaultpen(1);
  usepackage("amssymb");
  import geometry;

  // 
  point o = (8,0);
  pair ph = (10,0);
  pair pv = (5,8);

  pair pq3 = (3,1); //

  // 
  draw(o -- o+pv,Arrow(5bp,position=.4), Arrow(5bp,position=.5));
  draw(o+pv -- o+ph, Arrow(5bp,position=.6));
  draw(o+ph-pv-- o, Arrow(5bp,position=.6));
  draw(o+ph-pv-- o+ph, Arrow(5bp,position=.6));
  draw(o+ ph/2+(0,.2)--o+pv, Arrow(5bp,position=.4), Arrow(5bp,position=.5));
  draw(o+ph-pv -- o+ph/2+ (0,-.2));
  draw(o -- o+ph,Arrow(5bp,position=.3), Arrow(5bp,position=.4) );

   // 
  draw(o+pv -- o+ph+pq3, Arrow(5bp,position=.6));
  draw(o+ph-pv-- o+ph+pq3,Arrow(4bp,position=.4), Arrow(4bp,position=.5));
 draw(o+ph-- o+ph+pq3,Arrow(3bp,position=.3), Arrow(3bp,position=.4));

  // 
  //
   label("{\small  $q_1$}", o, 1*W);

  //
  label("{\small $q_2$}", o+ph, 2*dir(145));
   label("{\small $q_3$}", o+ph+pq3, .6*E);

   //
   label("{\small $p_1$}", o+pv, 1*N);

   //
  label("{\small $p_2$}", o+ph -pv, 1*S);

\end{asy}
\caption{A schematic view of the two tetrahedra glued along one face.}
\end{center}
\end{figure}

\subsection{The representation $\rho_2$}

In this section, we consider the second representation. Our main theorem is

\begin{thm}\label{maintheorem}
The representation $\rho_2$ is discrete and  is the holonomy of a branched spherical CR structure on the complement of the figure eight knot.
\end{thm}

The discreteness of the representation follows from the  observation that $\rho_2(\Gamma)$ is contained in a lattice. To prove the existence of a spherical CR structure on the complement of the figure eight knot,
it suffices to construct  two tetrahedra in the Heseinberg space with side pairings which allow the definition of a developing map.

The rest of this section is will be devoted to the construction of  the two tetrahedra in the Heseinberg space and to verify the conditions so that the developing map be well defined.  The difficulty of this construction is that we don't have a canonical way to
define the 2-skeleton.  The definitions of the faces are made so that they satisfy the necessary intersection properties.

We use half of a $\mathbb{C}$-circle to construct the segment between two given points.
For a given pair points in the Heisenberg space $p$ and $q$, we use the $[p,q]$ to denote the
segment connecting the two points with the direction from $p$ to $q$.

\subsection{The 0-skeleton and the side parings.}

The tetrahedra are $T_1:=[p_1,p_2,q_1,q_2]$ and $T_2:=[p_1,p_2,q_2,q_3]$,
where
$$
p_1=\infty, p_2=(0,0), q_1=(1,\sqrt{7}),
$$
$$
q_2=(\frac{5}{4}+i\frac{\sqrt{7}}{4},0), q_3=(\frac{1}{4}+i\frac{\sqrt{7}}{4},-\frac{\sqrt{7}}{2}).
$$

The side paring transformations are
$$
g_1: (q_2,q_1,p_1) \rightarrow (q_3,p_2,p_1)
$$
$$
g_2: (p_2,q_1, q_2) \rightarrow (p_1,q_2,q_3)
$$
$$
g_3: (q_1,p_2,p_1) \rightarrow (q_2,p_2,q_3)
$$

There are 6 tetrahedra around the edge $[p_2,p_1]$ (see Figure \ref{alongedge}) and $[p_2,q_2]$ respectively.  They are obtained by translating
$T_1$ and $T_2$.  They are:
$$
T_1, T_2, G_1(T_1), G_1G_3^{-1}(T_2), G_1G_3^{-1}G_2(T_1), G_1G_3^{-1}G_2G_1^{-1}(T_2),
$$
and respectively
$$
T_1,T_2,G_3(T_1),G_3G^{-1}_2(T_2) ,G_3G^{-1}_2(T_1), G_3G^{-1}_2G^{-1}_1(T_2).
$$

Following the side parings, it is easy to see the following:
\begin{enumerate}
  \item $G_1(q_2,q_1,p_1)=(q_3,p_2,p_1)$.
  \item $G_1G_3^{-1}(q_2,p_2,q_3)=(p_2,q_4,p_1)$ with
  $$q_4=G_1(p_2)=(-1,-\sqrt{7}).$$
  \item $G_1G_3^{-1}G_2(p_2,q_1,q_2)=(q_5,p_2,p_1)$ with
  $$q_5=G_1G_3^{-1}(p_1)=\left(\frac{-5}{4}+i\frac{\sqrt{7}}{4},0\right).$$
  \item $G_1G_3^{-1}G_2G_1^{-1}(q_3,p_2,p_1)=(p_1,p_2,q_6)$
  with
  $$q_6=G_1G_3^{-1}G_2(p_1)=\left(-\frac{1}{4}+i\frac{\sqrt{7}}{4},\frac{\sqrt{7}}{2}\right)$$
  and $G_1G_3^{-1}G_2G_1^{-1}(q_2)=q_1$.
  \item Since $G_1G_3^{-1}G_2G_1^{-1}G_3(q_1,p_2,p_1)=(q_1,p_2,p_1)$, $G_1G_3^{-1}G_2G_1^{-1}G_3=Id$.
\end{enumerate}
And,
\begin{enumerate}
  \item $G_3(p_2,q_1,q_2)=(p_2,q_2,p_3)$ with
 $$p_3=G_3(q_2)=\left(\frac{23}{32}+i\frac{5\sqrt{7}}{32},-\frac{\sqrt{7}}{16}\right).$$
  \item $G_3G^{-1}_2(p_1,p_2,q_2)=(p_2,p_4,q_2)$ with
$$p_4=G_3G^{-1}_2(p_2)=\left(\frac{5}{8}+i\frac{\sqrt{7}}{8},0\right).$$
  \item $G_3G^{-1}_2(p_1,q_1,q_2)=(p_2,p_5,q_2)$ with
$$p_5=G_3G^{-1}_2(q_1)=\left(\frac{3}{4}+i\frac{\sqrt{7}}{4},0\right).$$
  \item $G_3G^{-1}_2G^{-1}_1(p_1,q_2,q_3)=(p_2,q_1,q_2).$
  \item Since $G_3G^{-1}_2G^{-1}_1G_2(p_2,q_1,q_2)=(p_2,q_1,q_2)$, $G_3G^{-1}_2G^{-1}_1G_2=Id$.
\end{enumerate}

\subsection{ The 1-skeleton:}
\begin{figure*}[!t]
\centering
\subfigure[The space view]{\scalebox{0.4}{\includegraphics{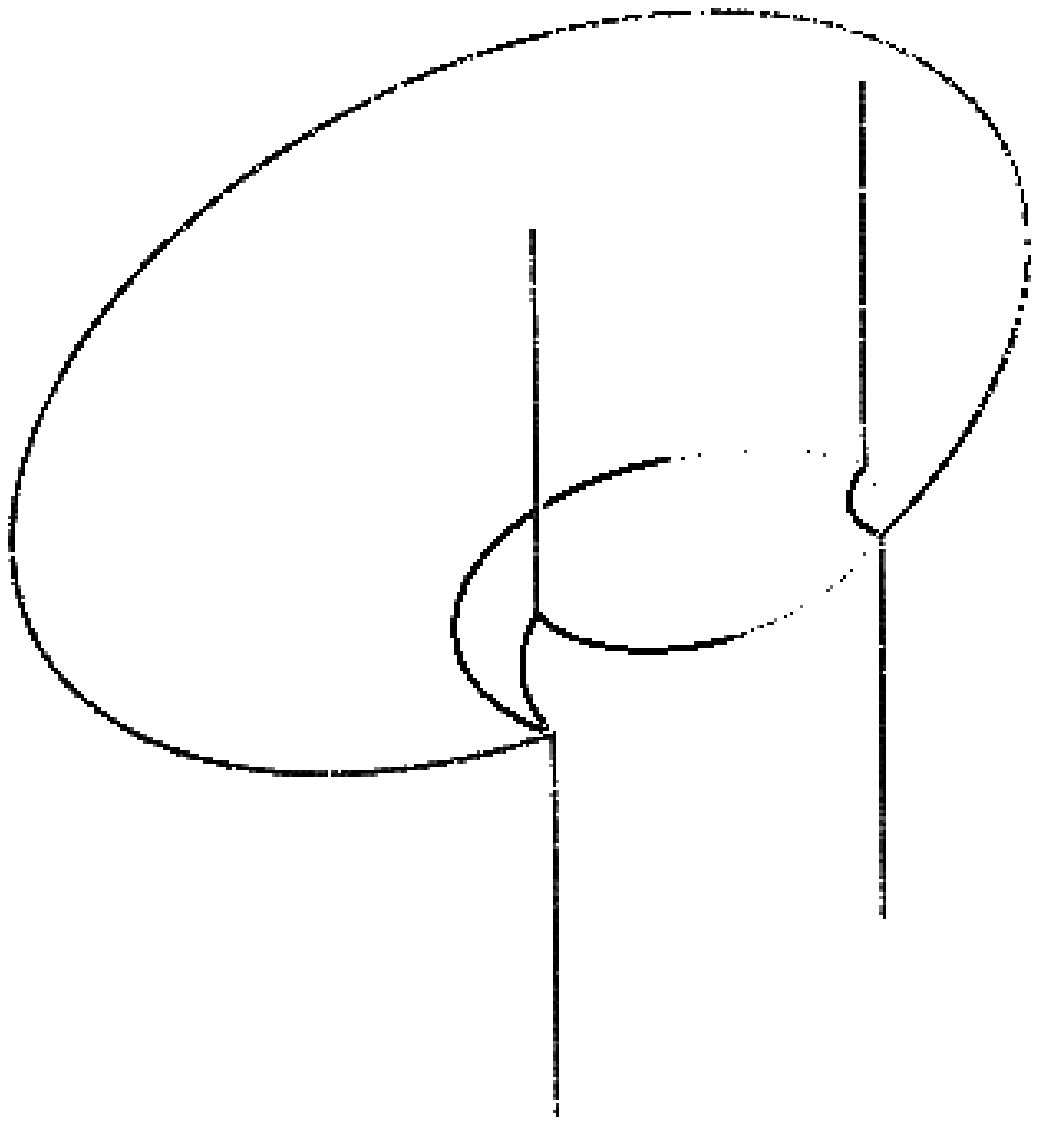}}}
\subfigure[The projected view]{\scalebox{0.2}{\includegraphics{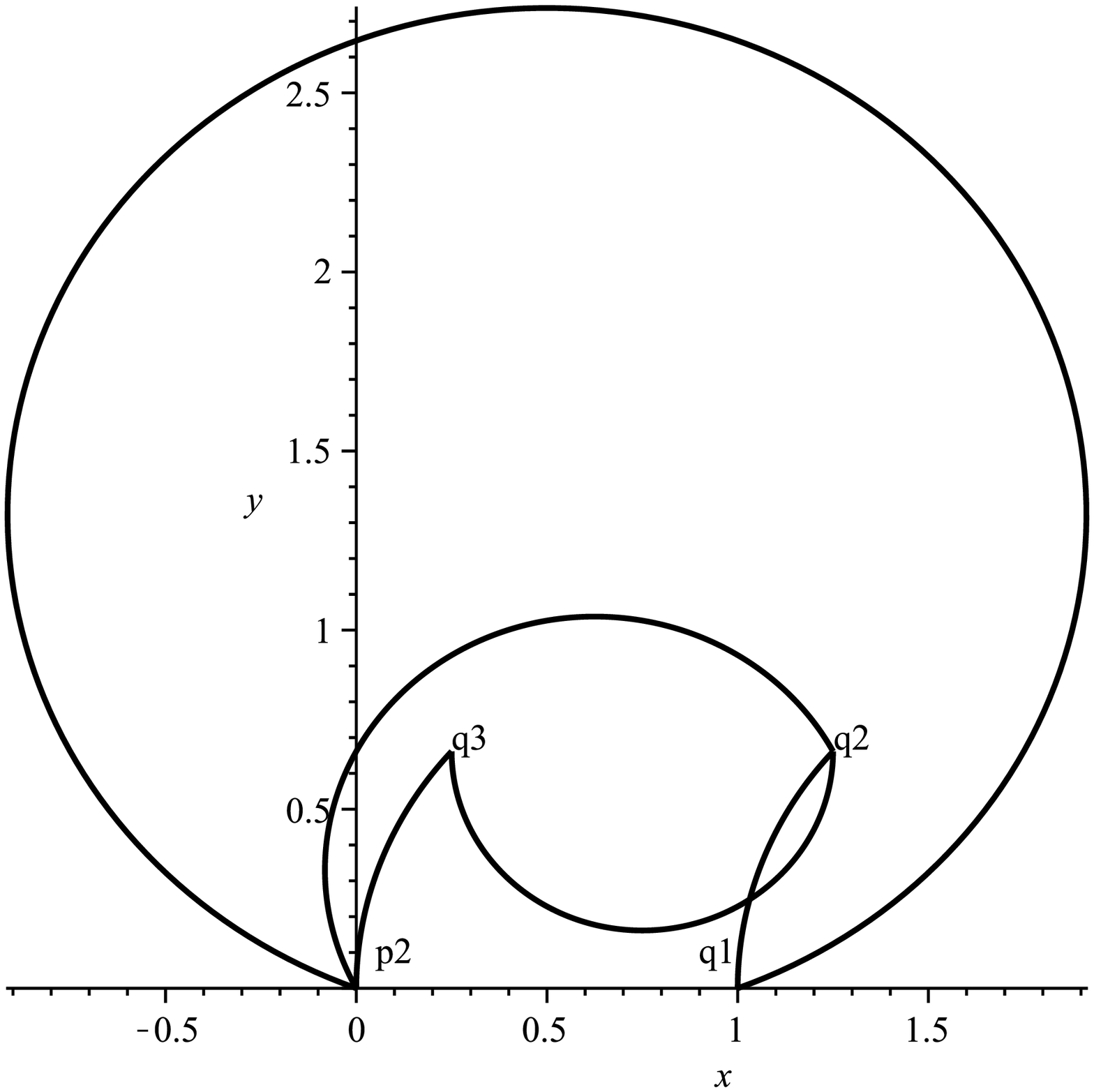}}}
\caption{The one skeleton of the tetrahedra}\label{f1}
\end{figure*}

In fact, considering the orientations of the edges, there are four possibilities for the choice of the one skeleton.
Here, we consider one choice given in Figure \ref{f1}. Precisely, $[p_2,p_1]=(0,t)$ with $t\leq 0$ and
$$[p_2,q_2]=\left(\frac{5+i\sqrt{7}}{8}+\frac{\sqrt{2}}{2}e^{i\theta},\frac{1}{8}\left(\sqrt{14}\cos(\theta)-5\sqrt{2}\sin(\theta)\right)\right),$$
where $\theta \in \left[\arccos\left(\frac{5\sqrt{2}}{8}\right),2\pi-\arccos\left(-\frac{5\sqrt{2}}{8}\right)\right]$. The other edges are determined from these two by applying
the side-pairings.

\subsection{The 2-skeleton}

In this subsection, we give the details of the construction of the faces of the two tetrahedra.

\subsubsection{Faces of $T_1$:}

\begin{figure}[h!!]
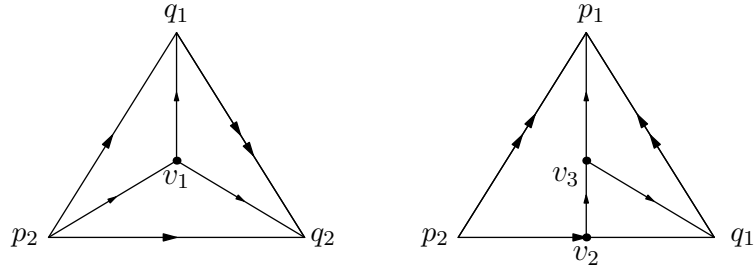
\label{faces1}
\begin{center}
\begin{asy}
  // 
  // 
  size(10cm);
  //defaultpen(1);
  usepackage("amssymb");
  import geometry;

  // 
  point o = (8,0);
  point oo = (24,0);
  pair ph = (10,0);
  pair pv = (5,8);
  point c= (13,3);
  point cc= (29,3);

  // 
  draw(o -- o+pv, Arrow(5bp,position=.5));
  draw(o+pv -- o+ph, Arrow(5bp,position=.6),Arrow(5bp,position=.5));
  draw(o -- o+ph,Arrow(5bp,position=.5) );
  dot("$v_1$",c,S);

path v1=c..o+pv;
path v2=c..o+ph;
path v3=o..c;
draw(v1,MidArrow(3bp));draw(v2,MidArrow(3bp));draw(v3,MidArrow(3bp));

//path uu2=A..u2..point(e1,.3);
//draw(uu,MidArrow(3bp));
//draw(uu2,MidArrow(3bp));

  // 
  label("{\small  $p_2$}", o, 1*W);
  label("{\small $q_2$}", o+ph, .6*E);
  label("{\small $q_1$}", o+pv, 1*N);

  // 
  draw(oo-- oo+pv,Arrow(5bp,position=.5),Arrow(5bp,position=.6));
  draw( oo+ph--oo+pv,Arrow(5bp,position=.5),Arrow(5bp,position=.6));
  draw(oo -- oo+ph,Arrow(5bp,position=.5));
dot("$v_3$",cc,SW);

path hor=oo -- oo+ph;
dot("$v_2$",point(hor,.5),S);
path hor=oo -- oo+ph;
//path u1=cc..u2..point(e1,.3);
path u1=cc..oo+pv;
path u2=cc..oo+ph;
path u3=point(hor,.5)..cc;
draw(u1,MidArrow(3bp));draw(u2,MidArrow(3bp));draw(u3,MidArrow(3bp));

   // 
  label("{ \small $p_2$}", oo, 1*W);
  label("{ \small $q_1$}", oo+ph, .6*E);
   label("{ \small $p_1$}", oo+pv, 1*N);

\end{asy}
\caption{A schematic view of faces $F(p_2,q_1,q_2)$ (left) and $F(p_1,p_2,q_1)$ (right).}
\end{center}
\end{figure}

\begin{figure}[h!!!]
\begin{center}
\begin{asy}
  // 
  // 
  size(4cm,0);
  //defaultpen(1);
 // usepackage("amssymb");
  //import geometry;
 real h=5;

path e2=(-h,0)..(-h,-5)..(0,-2*h);
path e1=(0,-2*h)..(h,-5)..(h,0);
path ie1=(h,0)..(h,-5)..(0,-2*h);//

draw(e2,MidArrow(5bp));
draw(ie1,Arrow(Relative(.7),size=1.7mm),Arrow(Relative(.6),size=1.7mm));

// On dŽfinit un rŽel et quatre points.
pair O=(0,0),A=(-h,0),B=(h,0),D=(0,-2*h);
// On trace le demi-cercle de diamtre [AB]
// en tournant dans le sens direct de B vers A.

draw(arc(O,B,A),MidArrow(5bp),Arrow(Relative(.6),size=1.8mm));
draw(arc(O,A,B),MidArrow(3bp));

// Et on ajoute les points.
dot("$p_1$",A,SW);
dot("$q_1$",B,SE);
dot("$q_2$",D,SW);

// Le 2-squelete
 pair u1=(-2,-h-.5);
pair u2=(-3.5,-h-1);

path uu=A..u1..point(e1,1);
path uu2=A..u2..point(e1,.3);
draw(uu,MidArrow(3bp));
draw(uu2,MidArrow(3bp));

draw(circle((-h+2*h/3,0),2*h/3),Arrow(3bp));
draw(circle((-h+h/2,0),h/2),Arrow(3bp));
draw(circle((-h+h/3,0),h/3),Arrow(3bp));

\end{asy}
\caption{A schematic view of the face $F(p_1,q_1,q_2)$.}
\end{center}
\end{figure}

\begin{figure}[h!!!]
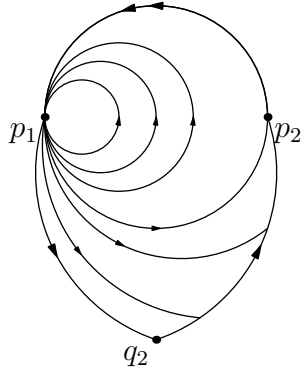

\begin{center}
\begin{asy}
  // 
  // 
  size(4cm,0);
  //defaultpen(1);
 // usepackage("amssymb");
  //import geometry;
 real h=5;

path e2=(-h,0)..(-h,-5)..(0,-2*h);
path e1=(0,-2*h)..(h,-5)..(h,0);

draw(e2,MidArrow(5bp));
draw(e1,MidArrow(5bp));

// On dŽfinit un rŽel et quatre points.
pair O=(0,0),A=(-h,0),B=(h,0),D=(0,-2*h);
// On trace le demi-cercle de diamtre [AB]
// en tournant dans le sens direct de B vers A.

draw(arc(O,B,A),MidArrow(5bp),Arrow(Relative(.6),size=1.8mm));
draw(arc(O,A,B),MidArrow(3bp));

// Et on ajoute les points.
dot("$p_1$",A,SW);
dot("$p_2$",B,SE);
dot("$q_2$",D,SW);

// Le 2-squelete
 pair u1=(-2,-h-.5);
pair u2=(-3.5,-h-1);

path uu=A..u1..point(e1,1);
path uu2=A..u2..point(e1,.3);
draw(uu,MidArrow(3bp));
draw(uu2,MidArrow(3bp));

draw(circle((-h+2*h/3,0),2*h/3),Arrow(3bp));
draw(circle((-h+h/2,0),h/2),Arrow(3bp));
draw(circle((-h+h/3,0),h/3),Arrow(3bp));

\end{asy}
\caption{A schematic view of the face $F(p_1,p_2,q_2)$.}
\end{center}
\end{figure}

We refer to Figure \ref{faces1} for a schematic description of the four faces.

\begin{enumerate}
  \item $F(p_2,q_1,q_2)$: choose $v_1=(\frac{3}{2}+i\frac{\sqrt{7}}{2},0)$ to be a center
   of the triangle $(p_2,q_1,q_2)$, then we define $F(p_2,q_1,q_2)$ to be the union of triangles $F(v_1,q_1,q_2)$, $F(p_2,q_1,v_1)$ and $F(p_2,q_2,v_1)$.
      \begin{itemize}
        \item $F(v_1,q_1,q_2)$ is the union of segments starting at $v_1$ and ending at the edge $[q_1,q_2]$;
        \item $F(p_2,q_1,v_1)$ is the union of segments starting at $p_2$ and ending at the edge $[v_1, q_1]$;
        \item $F(p_2,q_2,v_1)$ is the union of segments starting at $p_2$ and ending at the edge $[v_1, q_2]$.
      \end{itemize}
  \item $F(p_1,p_2,q_1)$: choose the point
  $$v_2=\left(\frac{1}{2}+i(\frac{\sqrt{7}}{2}+\sqrt{2}),-\sqrt{2}\right)\in [p_2,q_1,]$$
 and connect $v_2$ and $p_1$ by the edge $[v_2, p_1]$. Choose
  $$v_3=\left(\frac{1}{2}+i(\frac{\sqrt{7}}{2}+\sqrt{2}),-\sqrt{2}-6\sqrt{2}\right)\in [v_2,p_1],$$ then the face $F(p_1,p_2,q_1)$ is a union of faces $F(q_1,v_2,v_3)$, $F(p_1,q_1,v_3)$ and $F(p_1,p_2,v_2)$.
      \begin{itemize}
           \item $F(p_1,p_2,v_2)$ is the union of segments starting at each point of the segment $[p_2,v_2]$ and ending at $p_1$.
        \item $F(q_1,v_2,v_3)$ is the union of segments starting at each point of the segment $[v_2,v_3]$ and ending at $q_1$;
        \item $F(p_1,q_1,v_3)$ is the union of segments starting at each point of the segment $[v_3,q_1]$ and ending at $p_1$;
      \end{itemize}
  \item $F(p_1,q_1,q_2)$: It has two sub-faces, one is a triangle face which is the union of segments from $p_1$ to the edge $[q_1,q_2]$. The other one is a disc which is the union of $\mathbb{C}$-circles passing through $p_1$ and the half line $\{(1+it,\sqrt{7}): t\leq0\}$.
  \item $F(p_1,p_2,q_2)$: its construction is similar to the face $F(p_1,q_1,q_2)$. It also has two sub-faces, one is a union of segments from $p_1$ to the edge $[p_2,q_2]$, and the other is a disc which is the union of $\mathbb{C}$-circles passing through $p_1$ and the negative half of the $y$-axis in the Heisenberg space.
\end{enumerate}

\subsubsection{Faces of $T_2$:}
The faces of $T_2$ are all determined by the faces of $T_1$ by applying the side pairings.

\begin{enumerate}
  \item $F(p_1,q_2,q_3)$: Let
  $$v_4=G_2(v_1)=(\frac{3}{4}+i\frac{\sqrt{7}}{4},0).
  $$
  Since $F(p_1,q_2,q_3)=G_2(F(p_2,q_1,q_2))$, then $F(p_1,q_2,q_3)$ is a union of three faces, which are:
      \begin{itemize}
        \item $F(v_4,q_2,q_3)$ is the union of segments starting at $v_4$ and ending at the edge $[q_2,q_3]$;
        \item $F(p_1,v_4,q_2)$ is the union of the segments from  $p_1$ to the segment $[v_4,q_2]$;
         \item $F(p_1,v_4,q_3)$ is the union of the segments from $p_1$ to the segment $[v_4,q_3]$.
      \end{itemize}
  \item $F(q_3,p_2,q_2)$: Let
  $$v_5=G_3(v_2)\in [p_2,q_2].$$
  Connect $v_5$ and $q_3$ by the edge $[v_5,q_3]=G_3([v_2,p_1])$ and let $v_6=G_3(v_3)\in [v_5,q_3]$.  Since $F(q_3,p_2,q_2)=G_3(F(p_1,p_2,q_1))$), the face $F(q_3,p_2,q_2)$ is a union of three faces, which are:
      \begin{itemize}
        \item $F(q_3,p_2,v_5)$ is the union of the segments from the segment $[p_2,v_5]$ to $q_3$;
        \item $F(q_2,v_5,v_6)$ is the union of the segments from the segment $[v_5,v_6]$ to $q_2$.
        \item $F(q_3,v_6,q_2)$ is the union of the segments from the segment $[v_6,q_2]$ to $q_3$;
      \end{itemize}
  \item $F(p_1,p_2,q_2)$: It is the same as the definition of that face in the tetrahedron $T_1$.
  \item $F(p_1,p_2,q_3)$: From $F(p_1,p_2,q_3)=G_1(F(p_1,q_1,q_2))$, it is easy to see that the face $F(p_1,p_2,q_3)$ is the union of segments from
  $p_1$ to the edge $[p_2,q_3]$ and a disc which is a union of $\mathbb{C}-$circles passing through $p_1$ and the negative half of the y-axis.
\end{enumerate}

\subsection{The tetrahedra.}

In this subsection, we want to show that the faces of the tetrahedra constructed above define two tetrahedra.

Following the construction of the 2-skeleton, it is easy to show each face is embedded.
\begin{lem}
Each face of the two tetrahedra defined in the above section is topologically a disc in the Heisenberg space.
\end{lem}



\begin{lem}\label{L2}
The tetrahedron $T_1$ defined above is homeomorphic to a tetrahedron.
\end{lem}

\begin{lem}\label{L1}
The tetrahedron $T_2$ defined above is homeomorphic to a generalized tetrahedron.
\end{lem}

\begin{lem}\label{L3}  $T_1 \bigcap T_2=F(p_1,p_2,q_2)$.
\end{lem}

From the definition of $T_1$ and $T_2$, and the above lemmas, we have
\begin{lem}
$G_1,G_2,G_3$ are side parings of the union $T_1\bigcup T_2$.
\end{lem}

\begin{pro}
The quotient space of $T_1\bigcup T_2 - \{vertices\}$ under the side parings $G_1,G_2,G_3$ is the complement of the figure eight knot.
\end{pro}

\subsection{The structure around the edges}
The quotient of $T_1\bigcup T_2 $ by the side parings has two edges, represented by $[p_2,p_1]$ and $[p_2,q_2]$.
The purpose of this subsection is to show that the neighborhood around those edges covers a neighborhood
of half of the $t$-axis in the Heisenberg space. The phenomenon is similar as that in the subsection 6.4 of \cite{fal08}.

\subsubsection{The neighborhood around $[p_2,p_1]$.}

We know that the neighborhood around $[p_2,p_1]$ is a union of the neighborhoods
contained in
$$T_1=[p_1,p_2,q_1,q_2],$$
$$T_2=[p_1,p_2,q_2,q_3],$$
$$T_3=G_1(T_1)=[p_1,q_4,p_2,q_3],$$
$$T_4=G_1G_3^{-1}(T_2)=[q_5,q_4,p_2,p_1],$$
$$T_5=G_1G_3^{-1}G_2(T_1)=[q_6,q_5,p_2,p_1],$$
and
$$T_6=G_1G_3^{-1}G_2G_1^{-1}(T_2)=[q_6,p_2,q_1,p_1].$$

From the above six tetrahedra, we know that the six faces with the same edge $[p_2,p_1]$ are $F(p_2,p_1,q_j)$, where $j=1..6$.
By arguing as in \cite{fal08}, it is easy to see that each pair of consecutive tetrahedra $T_j$ and $T_{j+1}$  match monotonically along the matching face $F(p_2,p_1,q_{j+1})$.

\begin{figure}[h!!!]
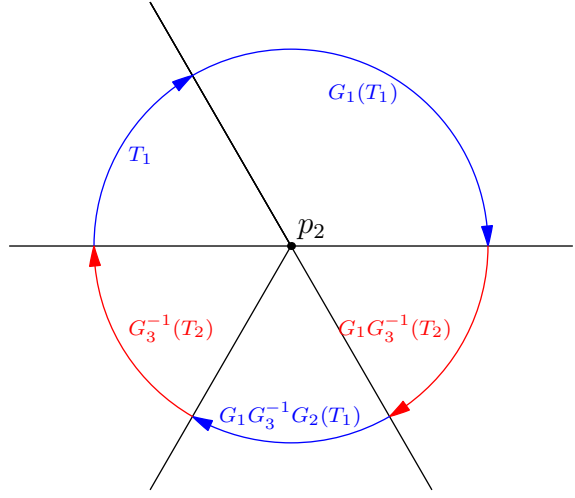

\begin{center}
\begin{asy}

 size(7.5cm,0);
import geometry;

pair O=(0,0);
dot("$p_2$",O,NE);
real x=0.7;
//draw(unitcircle);

real a1=180, a2=120,a3=120,a4=0,a5=-60,a6=-120;
pair q1=dir(a1),  q2=dir(a2),  q3=dir(a3), q4=dir(a4), q5=dir(a5), q6=dir(a6);
pair q4e=dir(-.01);

draw(O--q1);draw(O--q2);draw(O--q2);draw(O--q3);draw(O--q4);draw(O--q5);draw(O--q6);

draw("$^{T_1}$",arc(q1,O,q2,x),blue,Arrow);
draw("$^{G_1(T_1)}$",arc(q3,O,q4,x),blue,Arrow);
draw("$^{G_1G_3^{-1}(T_2)}$",arc(q4e,O,q5,x),red,Arrow);
draw("$^{G_1G_3^{-1}G_2(T_1)}$",arc(q5,O,q6,x),blue,Arrow);
draw("$^{G_3^{-1}(T_2)}$",arc(q6,O,q1,x),red,Arrow);

\end{asy}
  \caption{A schematic picture of a  neighborhood around the edge $[p_1,p_2]$, where the segments stand for the
  faces with the common edge $[p_1,p_2]$ denoted by the common intersection point $p_2$, and the arcs and the regions  between two segments stand for the neighborhoods
  contained in one tetrahedron. Here $T_2$ degenerates to a subface of $F(p_2,p_1,q_2)$, which is the same subface of $F(p_2,p_1,q_3)$.}\label{f18}
\end{center}
\end{figure}

Let $N_j,j=1..6$ denote the neighborhoods around the edge $[p_2,p_1]$ contained in the six tetrahedra $T_j$.
By analyzing the positions of those neighborhoods in the Heisenberg space (see the Figure \ref{f18} for a schematic description)
we have the following proposition:
\begin{pro}
The union $\bigcup N_j$ forms a standard tubular neighborhood of $[p_2,p_1]$ in the Heisenberg space.
\end{pro}

 \subsubsection{The neighborhood around $[p_2,q_2]$.}

The neighborhood around $[p_2,q_2]$ is a union of the neighborhoods contained in the six tetrahedra
$$T'_1=T_1=[p_1,p_2,q_1,q_2],$$
$$T'_2=T_2=[p_1,p_2,q_2,q_3],$$
$$T'_3=G_3(T_1)=[q_3,p_2,q_2,p_3],$$
$$T'_4=G_3G^{-1}_2(T_2)=[p_2,p_4,q_2,p_3],$$
$$T'_5=G_3G^{-1}_2(T_1)=[p_2,p_4,p_5,q_2]$$
and
$$T'_6=G_3G^{-1}_2G^{-1}_1(T_2)=G^{-1}_2(T_2)=[p_2,p_5,q_1,q_2].$$

\begin{figure}[h!!!]
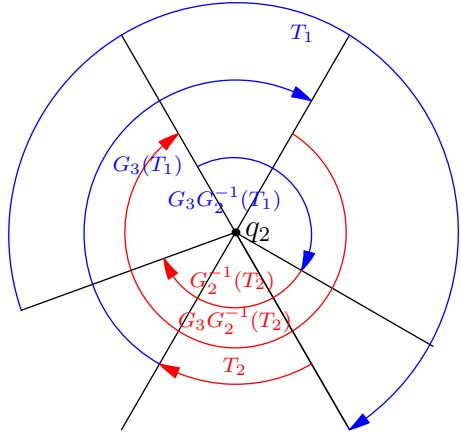

\begin{center}
\begin{asy}

 size(6cm,0);
import geometry;

pair O=(0,0);
dot("$q_2$",O,E);
real x=2/3;
//draw(unitcircle);

real a1=200, a2=300,a3=-120,a4=60,a5=120,a6=-30;
pair q1=dir(a1),  q2=dir(a2),  q3=dir(a3), q4=dir(a4), q5=dir(a5), q6=dir(a6);

draw(O--q1);draw(O--q2);draw(O--q2);draw(O--q3);draw(O--q4);draw(O--q5);draw(O--q6);

draw("$^{T_1}$",q1..q5..q4..q2,blue,Arrow);
draw("$_{T_2}$",arc(q2,O,q3,x),red,Arrow);
draw("$^{G_3(T_1)}$",x*q3..x*q1..x*q5..x*q4,blue,Arrow);
draw("$^{G_3G_2^{-1}(T_2)}$",1/2*q4..1/2*q2..1/2*q3..1/2*q5,red,Arrow);
draw("$^{G_3G_2^{-1}(T_1)}$",1/3*q5..(1/3,0)..1/3*q6,blue,Arrow);
draw("$^{G_2^{-1}(T_2)}$",1/3*q6..1/3*q3..1/3*q1,red,Arrow);

\end{asy}
 \caption{A schematic picture showing the neighborhoods corresponding to each tetrahedron. The segments stand for the faces with the common  edge $[p_2,q_2]$ ( which is represented by the point $q_2$) and the arcs and the region between two segments denote the neighborhoods contained in one tetrahedron.}\label{f24}\end{center}
\end{figure}

Let $N'_j$, $j=1..6$ denote the neighborhood around the edge $[p_2,q_2]$ contained in the corresponding tetrahedron $T'_j$.
One can analysis the positions of those tetrahedra. In fact, it suffice to analysis the faces containing the same edge $[p_2,q_2]$ and the intersections
with a tubular neighborhood of $[p_2,q_2]$. (See the Figure \ref{f24} for an abstract description of the position of the neighborhoods). By a similar argument in \cite{fal08}, we can conclude
\begin{pro}
The union $\bigcup N'_j$ forms a neighborhood covering three
times a standard tubular neighborhood of the edge $[p_2,q_2]$ in the Heisenberg space.
\end{pro}

\begin{remk}
We correct a statement in \cite{fal08}.  In fact  the union of the neighborhoods contained in the tetrahedra around the edge $[p_2,p_4]$  (in the case of the first representation discussed there) forms a standard neighborhood of this edge, and not a three times cover as announced in the paper.
\end{remk}

\section{Appendix}

In this section, we give the proofs of Lemma \ref{L2}, Lemma \ref{L1} and Lemma \ref{L3} contained in the subsection 6.5.

Let $\Pi: \mathfrak{N} \rightarrow \mathbb{C}$ be vertical projection map from the Heisenberg space onto the z-plane.  When describing projections in this section,
we will use the same notation for a point in the Heisenberg group and its projection in the $z$-plane.

\subsection{Proof of  Lemma \ref{L2}:}

It suffice to show that each pair of  faces only intersect at their common edge. It is well known that any $\mathbb{C}-$circle passing through
the point at infinity is a vertical line in the Heisenberg space. Hence any segment with $p_1$ as an endpoint will project to a point on the z-plane.

\begin{figure}
  \centering
{\scalebox{0.3}{\includegraphics{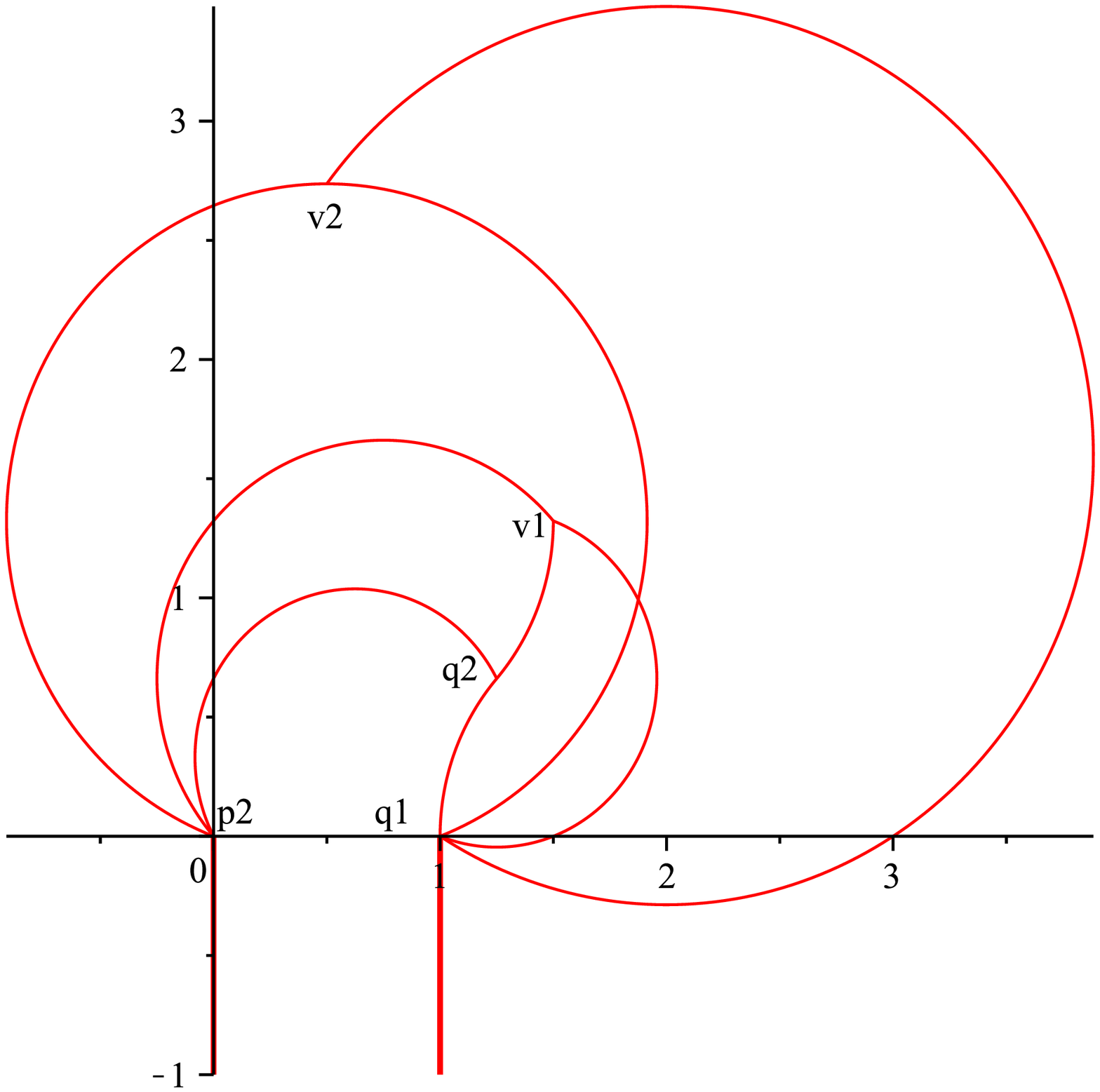}}}
  \caption{The projections of the faces of $T_1$
  }
  \label{f11}
\end{figure}

First, we analyze the projections on the z-plane of  the projections of the faces of the tetrahedron $T_1$.
It is easy to determine their projections (see the Figure \ref{f11} ):
\begin{itemize}
  \item $\Pi(F(p_1,p_2,q_2))$: The union of the (circle) curve $p_2 q_2$ and the negative $y$-axis starting at $p_2$ ;
  \item $\Pi(F(p_1,q_1,q_2))$: The union of the (circle) curve $q_1 q_2$ and the half-line parallel to the $y$-axis starting at $q_1$ ;
  \item $\Pi(F(p_1,p_2,q_1))$: The union of $p_2 v_2$ and the region between the two curves connecting $v_2$ and $q_1$;
  \item $\Pi(F(p_2,q_1,q_2))$: It is the union of the triangles $(p_2,v_1,q_2)$ and $(v_1,q_1,q_2)$ and the  curves from the point $p_2$ to the curve $v_1 q_1$.
\end{itemize}

\begin{figure}
  \centering
{\scalebox{0.3}{\includegraphics{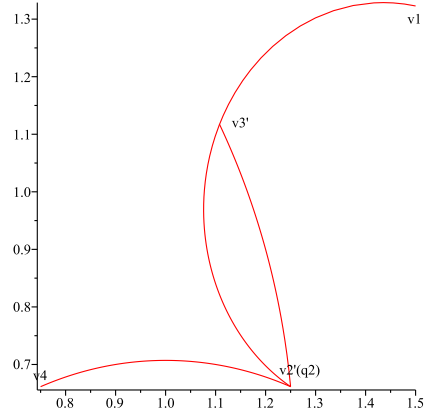}}}
  \caption{The projections $\Pi(F(p_1,q_2,v_4))$ is the curve $v_4 q_2$, and $\Pi(F(q_2,v'_2,v'_3))$ is the region bounded by
  the two curves connecting $v'_2$ and $v'_3$.}\label{f12}
\end{figure}

 The only one we have to
check carefully is
$$F(p_2,q_1,q_2) \bigcap F(p_1,p_2,q_1)=[p_2,q_1],$$
since the others obviously intersect at their common edge. Recall that both of the faces contain three sub-faces, so it suffices to prove
$$F(p_2,q_1,v_1) \bigcap F(q_1,v_2,v_3)=[v_2,q_1]$$
 since $v_2\in [p_2,q_1]$. As it is not easy to see this from their projections, we consider the images of these two faces by the transformation $G_2$ which will transform the point $p_2$ to the point at infinity $p_1$.
\begin{equation*}
\begin{split}
G_2(F(p_2,q_1,v_1) \bigcap F(q_1,v_2,v_3))
&=G_2(F(p_2,q_1,v_1)) \bigcap G_2(F(q_1,v_2,v_3))\\
&=F(p_1,q_2,v_4)\bigcap F(q_2,v'_2,v'_3),
\end{split}
\end{equation*}
where
$$v'_2=G_2(v_2)=\left(\frac{5}{4}+i\frac{\sqrt{7}}{4},\sqrt{2} \right)\in G_2([p_2,q_1])=[p_1,q_2]$$
 and
 $$v'_3=G_2(v_3)=\left( \frac{40+\sqrt{14}}{32+2\sqrt{14}}+i\frac{5\sqrt{2}+14\sqrt{7}}{32+2\sqrt{14}},-\frac{\sqrt{2}+2\sqrt{7}}{32+2\sqrt{14}} \right)\in G_2([v_2,p_1])=[v'_2,v_1].$$
It can be seen that
$$F(p_1,q_2,v_4)\bigcap F(q_2,v'_2,v'_3)=[v'_2,q_2]=G_2([v_2,q_1])$$
 by analyzing their projections (see the Figure \ref{f12}), which completes our proof.

\subsection{Proof of  Lemma \ref{L1}:}

\begin{figure}
  \centering
{\scalebox{0.3}{\includegraphics{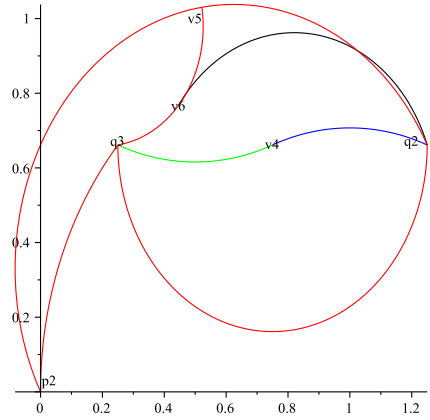}}}
  \caption{The projections of the faces of $T_2$
  }\label{f13}
\end{figure}

As in the proof of the above lemma, we first consider the projections of the faces given in Figure \ref{f13}:
\begin{itemize}
  \item $\Pi(F(p_1,p_2,q_2))$: It is the union of the circle segment $p_2 q_2$ and  the negative half y-axis which is the projection of the disc part;
  \item $\Pi(F(p_1,p_2,q_3))$: It is the union of the circle segment $p_2 q_3$ and  the negative half y-axis;
  \item $\Pi(F(p_1,q_2,q_3))$: It is the triangle $(v_4,q_2,q_3)$;
  \item $\Pi(F(p_2,q_2,q_3))$: It is the union of  $\Pi(F(p_2,v_5,q_3))$,  $\Pi(F(v_6,q_2,q_3))$ and  $\Pi(F(q_2,v_5,v_6))$ which is the union of (circle) curves from the point $q_2$ to the (circle) curve $v_5 v_6$.
\end{itemize}
The intersections of each pair of  faces are easily obtained
except
$$F(p_1,q_2,q_3)\bigcap F(p_2,q_2,q_3)=[q_2,q_3]$$
and
 $$F(p_2,q_2,q_3)\bigcap F(p_2,q_2,p_1)=[p_2,q_2].$$

\begin{figure}
  \centering
{\scalebox{0.3}{\includegraphics{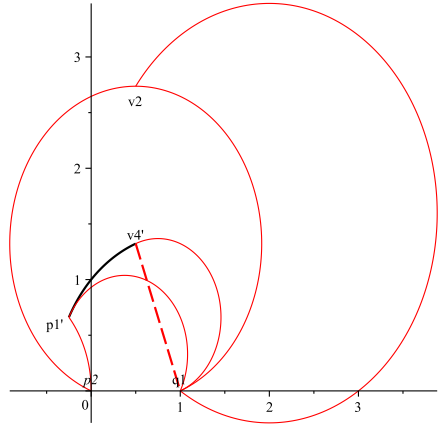}}}
  \caption{The projections of the faces $F(p_1',q_1,p_1)$, $ F(p_2,q_1,p_1)$ and $F(p_2,q_1,p_1')$.}\label{f14}
\end{figure}

The first one can be obtained by considering their images under $G_3^{-1}$. We have to show
\begin{equation*}
\begin{split}
G_3^{-1}(F(p_1,q_2,q_3)\bigcap F(p_2,q_2,q_3))
&=G_3^{-1}(F(p_1,q_2,q_3)) \bigcap G_3^{-1}(F(p_2,q_2,q_3))\\
&=F(G_3^{-1}(p_1),q_1,p_1)\bigcap F(p_2,q_1,p_1)\\
&=[q_1,p_1]=G_3^{-1}([q_2,q_3]).
\end{split}
\end{equation*}
Let $$p_1'=G_3^{-1}(p_1)=\left(-\frac{1}{4}+i\frac{\sqrt{7}}{4},\frac{\sqrt{7}}{2}\right)$$
 and $$v_4'=G_3^{-1}(v_4)=\left(\frac{1}{2}+i\frac{\sqrt{7}}{2},0\right).$$
It suffice to show that
$$ F(v_4',q_1,p_1)\bigcap F(q_1,v_2,p_1)=[q_1,p_1],$$
 since the other two sub-faces $F(p_1,v_4,q_2)$ and $F(p_1,v_4,q_3)$ of $F(p_1,q_2,q_3)$ do not intersect the face $F(p_2,q_2,q_3)$. This can be verified by analyzing their projections in Figure \ref{f14}, where
the projections of $F(v_4',q_1,p_1)$ lie in the region between the straight line and the circle segment with the same endpoints $v_4'$ and $q_1$.

To prove
$$F(p_2,q_2,q_3)\bigcap F(p_2,q_2,p_1)=[p_2,q_2],$$
 is equivalent to show
\begin{equation*}
\begin{split}
G_3^{-1}(F(p_2,q_2,q_3))\bigcap G_3^{-1}(F(p_2,q_2,p_1))
&=F(p_2,q_1,p_1)\bigcap F(p_2,q_1,p_1')\\
&=[p_2,q_1].
\end{split}
\end{equation*}
The result is clear analysing the projections in Figure \ref{f14}, since the projection of $F(p_2,q_1,p_1')$ lies in the triangle $(p_1',p_2,q_1)$.

At last, we have to mention that the intersection of $F(p_1,p_2,q_2)$ and $F(p_1,p_2,q_3)$ is a disc, not only an edge ($T_2$ is a generalized tetrahedron).

\subsection{Proof of Lemma \ref{L3}:}

\begin{figure}[!]
  \centering
{\scalebox{0.3}{\includegraphics{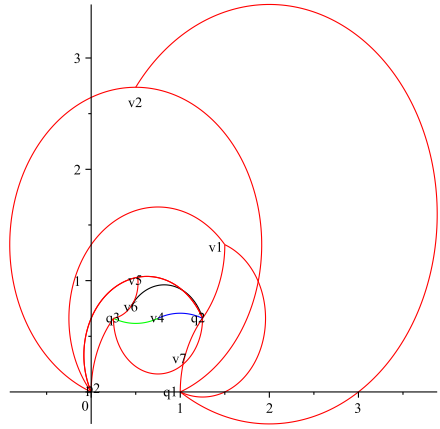}}}
  \caption{The projection of the subfaces of $T_1$ and $T_2$}\label{f15}
\end{figure}

\begin{figure}[!]
  \centering
{\scalebox{0.3}{\includegraphics{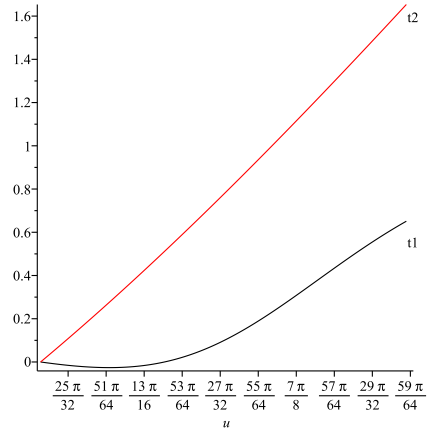}}}
  \caption{The height comparison, where $t_1,t_2$ denote the height in $[v_7,q_2]$ and $F(v_4,q_2,q_3)$ respectively.}\label{f16}
\end{figure}

According to the projections of the faces of the two tetrahedra given in Figure \ref{f15}, we only need to prove the following cases in detail.
\begin{itemize}
  \item $F(p_1,q_1,q_2)\bigcap F(p_1,q_2,q_3)=[p_1,q_2]$.
  Recall that the face $F(p_1,q_2,q_3)$ has three parts. It is easy to see
  $$F(p_1,q_1,q_2)\bigcap F(p_1,q_2,v_4)=[p_1,q_2]$$
   and
   $$F(p_1,q_1,q_2)\bigcap F(p_1,v_4,q_3)=p_1.$$
    Therefore, it is suffice to show
  $$
  F(v_4,q_2,q_3)\bigcap F(q_1,q_2,p_1)=q_2,
  $$
  which follows by comparing the height functions of the two faces.
  From their projections, we only need to compare the height of the parts where they have
  intersected projections, i.e. the segment $[v_7,q_2]\subset [q_1,q_2]$
  where $v_7=\left(\frac{33}{32}+i\frac{3 \sqrt{7}}{32},\frac{5\sqrt{7}}{8}\right)$. More precisely, write the x and y coordinates of
  $$[v_7,q_2]=\left(2+e^{i\theta},\sqrt{7}-4\sin(\theta)\right)$$
  where
  $$\theta \in \left[\pi-\arcsin\left(\frac{\sqrt{7}}{4}\right),\pi-\arcsin\left(\frac{3\sqrt{7}}{32}\right)\right]$$
  into the parametrization of the face
  $$
  F(v_4,q_2,q_3)=\left\{
                   \begin{array}{ll}
                     (x-x_0)^2+(y-y_0)^2=1/2 \\
                      t=t_0+2(y_0 x-x_0 y)
                   \end{array}
                 \right.
  $$
  where
  $$ \left\{
       \begin{array}{ll}
          x_0=\left(\cos(\varphi)+\sqrt{7}\sin(\varphi)+3\right)/4 \\
          y_0=\left(-\sqrt{7}\cos(\varphi)+\sin(\varphi)+\sqrt{7}\right)/4 \\
          t_0=\left(\sqrt{7}\cos(\varphi)+\sin(\varphi)\right)/2
       \end{array}
     \right.
  $$
  with $\varphi \in [\pi, 2\pi]$,
  we can get the height function $t_1=t_1(\theta)$ as a function of $\theta$.
  Let $t_2=\sqrt{7}-4\sin(\theta)$, then we can compare these two height functions (see Figure \ref{f16})
  so that the height of $[v_7,q_2]$ is bigger than that in $F(v_4,q_2,q_3)$.
  This implies that $F(v_4,q_2,q_3)$ and $F(q_1,q_2,p_1)$ only intersect at the point $q_2$.
  \item $F(p_2,q_1,q_2)\bigcap F(p_2,q_2,q_3)=[p_2,q_2]$.
  \item $F(p_2,q_1,q_2)\bigcap F(p_1,q_2,q_3)=q_2$.
\end{itemize}

\begin{figure}[!]
  \centering
{\scalebox{0.3}{\includegraphics{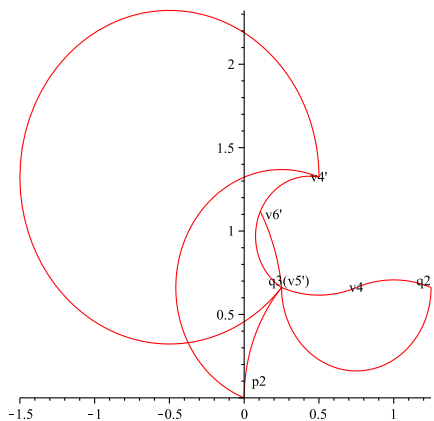}}}
  \caption{Projections of $F(p_1,q_2,q_3)$, $F(p_2,q_3,v_4')$ and $F(q_3,v_5',v_6').$
  }\label{f17}
\end{figure}

The last two can be proved by a similar argument as in the proof of Lemma \ref{L2} and Lemma \ref{L1}.
More precisely, we consider their images under the action of $G_2$. Recall that
$$G_2(F(p_2,q_1,q_2))=F(p_1,q_2,q_3)$$
 and  $v_1=G_2(p_1)$.
  Let $v'_4=G_2(q_3)$, then
 $$G_2(F(p_1,q_2,q_3))=F(v_1,q_3,v_4')$$ and
 $$G_2(F(p_2,q_2,q_3))=F(p_1,q_3,v_4').$$
 Recall that each of these faces has three parts, according to the projected view in Figure \ref{f15} we only need to check the intersections of their subfaces
 $$F(p_2,v_1,q_2)\bigcap F(v_5,v_6,q_2)=[v_5,q_2]$$
  and
$$F(v_1,q_1,q_2)\bigcap F(v_4,q_2,q_3)=q_2.$$
These can be proved by analyzing the projections of
their images by $G_2$ in Figure \ref{f17}.
Precisely, let
$$
v_5'=G_2(v_5)=\left(\frac{1}{4}+i\frac{\sqrt{7}}{4},\frac{\sqrt{2}}{8+2\sqrt{14}} \right),
$$
$$
v_6'=G_2(v_6)=\left( \frac{8-\sqrt{14}}{32+2\sqrt{14}}+i\frac{5\sqrt{2}+14\sqrt{7}}{32+2\sqrt{14}},-\frac{5\sqrt{2}+6\sqrt{7}}{32+2\sqrt{14}} \right)
$$
and recall that $p_2=G_2(v_4)$. Then their projections are:
\begin{itemize}
  \item $\Pi(F(p_1,q_2,q_3))$ is the triangle $(v_4,q_2,q_3)$;
  \item $\Pi(F(p_2,q_3,v_4'))$ is the union of (circle) curves from  the curve $p_2 q_3$ to the point $v'_4$.  This projection is more complicated but lies outside the triangle $(v_4,q_2,q_3)$;
  \item $\Pi(F(q_3,v_5',v_6'))$ is the region bounded by the two curves connecting the points $v'_5$ and $v'_6$.
\end{itemize}
Observe that the points $v'_5$ and $q_3$ denote the same points on the z-plane.

\end{document}